\documentclass{SCAEOL}
\numberwithin{equation}{section}

\usepackage{algorithm,algorithmic}
\usepackage{graphicx,subfig}
\usepackage{color}
\usepackage{csquotes}

\usepackage{amsopn}
\DeclareMathOperator{\diag}{diag}
\DeclareMathOperator*{\argmin}{argmin}

\newcommand{\R}{\mathbb{R}}
\newcommand{\T}{\mathcal{T}}
\newcommand{\vecq}{\mathbf{q}}

\newcommand{\vecz}{z}
\newcommand{\vecw}{w}
\newcommand{\vecmu}{\boldsymbol{\lambda}}

\usepackage[capitalize,nameinlink]{cleveref}
\crefname{figure}{Figure}{Figures}
\crefname{table}{Table}{Tables}
\crefname{proposition}{Proposition}{Propositions}
\crefname{lemma}{Lemma}{Lemmas}
\crefname{corollary}{Corollary}{Corollaries}
\crefname{definition}{Definition}{Definitions}

\begin{document}

\Year{2013} %
\Month{January}
\Vol{56} %
\No{1} %
\BeginPage{1} %
\EndPage{XX} %
\AuthorMark{First1 L N {\it et al.}}

\title[Truncated Regularization Framework]{A General Truncated Regularization Framework for Contrast-Preserving Variational Signal and Image Restoration: Motivation and Implementation}{}


\author[]{WU ChunLin}{}
\author[]{LIU ZhiFang}{Corresponding author}
\author[]{WEN Shuang}{}

\address[]{School of Mathematical Sciences, Nankai University, Tianjin {\rm 300071}, China}
\Emails{wucl@nankai.edu.cn,
liuzhifang0628@gmail.com, wenshuang13@126.com}

\maketitle


 {\begin{center}
\parbox{14.5cm}{\begin{abstract}
Variational methods have become an important kind of methods in signal and image restoration - a typical inverse problem. One important minimization model consists of the squared $\ell_2$ data fidelity (corresponding to Gaussian noise) and a regularization term constructed by a potential function composed of first order difference operators. It is well known that total variation (TV) regularization, although achieved great successes, suffers from a contrast reduction effect. Using a typical signal, we show that, actually all convex regularizers and most nonconvex regularizers have this effect. With this motivation, we present a general truncated regularization framework. The potential function is a truncation of existing nonsmooth potential functions and thus flat on $(\tau,+\infty)$ for some positive $\tau$. Some analysis in 1D theoretically demonstrate the good contrast-preserving ability of the framework. We also give optimization algorithms with convergence verification in 2D, where global minimizers of each subproblem (either convex or nonconvenx) are calculated. Experiments numerically show the advantages of the framework.
 \vspace{-3mm}
\end{abstract}}\end{center}}



\renewcommand{\baselinestretch}{1.2}
\begin{center} \renewcommand{\arraystretch}{1.5}
{\begin{tabular}{lp{0.8\textwidth}} \hline \scriptsize
{\bf Citation:}\!\!\!\!&\scriptsize First1 L N, First2 L N, First3 L N.  SCIENCE CHINA Mathematics  journal sample. Sci China Math, 2013, 56, doi: 10.1007/s11425-000-0000-0\vspace{1mm}
\\
\hline
\end{tabular}}\end{center}

\baselineskip 11pt\parindent=10.8pt  \wuhao

\section{Introduction}
A digital signal (1D) or image (2D) records the intensity information of a scene or object in our real world. Mathematically it can be regarded as a uniform sampling from a univariate or multivariate function. Due to some limitations of the sampling procedure, the signal or image data usually contain degradations such as noise and blur. Signal and image restoration is to recover the underlying clean and clear data, which is a fundamental task in image data processing. The most important and usually difficult target is to keep the edges (discontinuities) and their contrasts.

In recent decades variational methods (energy minimization methods) have become an important kind of methods for this typical inverse problem \cite{Tikhonov1977Solutions}, especially when the degradation operator is not an identity. The energy function usually contains two terms: regularization term and data fidelity term (see Section 2 for details). In a very common and important case with Gaussian noise, which will be considered in this paper, the data fidelity term is a squared $\ell_2$ norm. A basic choice for the regularization term is constructed by a regularization function composed of the discrete gradient operator, i.e., the first order finite difference operator. (Of course one may consider the second order, third order finite differences or wavelet frames.) The key of the method is therefore to choose a suitable regularization function. So far, several regularization functions have been proposed, including convex ones, smooth ones, and nonconvex and nonsmooth ones, in either the variational or statistical setting.

The convex or smooth potential functions were proposed in \cite{Rudin1992Nonlinear,Bouman1993generalized,Green1990Bayesian} or  \cite{Geman1987Statistical,Hebert1989generalized,Leclerc1989Constructing,Perona1990Scale,Vogel1996Iterative},\textit{etc}. For a smooth regularization based method, the Euler-Lagrange equation is a good way to analyze the variational model; see, e.g., \cite{Aubert2006Mathematical} and \cite{Charbonnier1997Deterministic} where some useful conditions on the regularization function were suggested to pursuit edge preservation. It is well known that the both convex and smooth $\ell_2$ regularization over-smoothes edges. A typical nonsmooth but convex potential function is the $\ell_1$ norm, which, when composed of the gradient operator, is the popular total variation (TV) regularization. Due to its edge preservation ability, TV attracted much attention and has many extensions; see the references in \cite{Scherzer2015Handbook}. Convex variational models often benefit from efficient algorithms~\cite{Wen2005new,Goldstein2009Split,Wu2010Augmented,Wen2016Primal}. However, as Meyer's remarkable example \cite{Meyer2002Oscillating} shows, TV regularization suffers from a contrast reduction effect.

Earlier nonconvex and nonsmooth regularization methods have been presented in \cite{Geman1984Stochastic,Blake1987Visual,Leclerc1989Constructing,Geman1992Constrained,Geman1995Nonlinear,Mumford1989Optimal,Saquib1998ML}, \textit{etc}. In these methods, some auxiliary variables were used to serve as markers of edges or some nonconvex and nonsmooth regularization functions were used to produce neat edges. In the latter case, a series of papers \cite{Nikolova2004Variational,Nikolova2005Analysis,Durand2006Stability,Durand2006Stabilitya,Nikolov2008Analytical,Nikolov2013Description} systematically  analyzed the properties of the minimizers of nonconvex and nonsmooth variational models. As noticed before that model validation is difficult in general, these analysis shows the advantages of nonconvex and nonsmooth models in image restoration, among which the most important one is the capability to produce neat edges. Very recently nonconvex and nonsmooth regularization obtained much attention in variational image restoration~\cite{Nikolova2008Efficient,Nikolova2010Fast,Fornasier2010Iterative,Xu2012L12,Dong2013Efficient,Lou2015Weighted,Hintermueller2013Nonconvex,Ochs2014iPiano,Bian2015Linearly,Kolmogorov2016Total,Lanza2016Constrained,Ochs2015Iteratively}.

In this paper we study variational restoration models from another viewpoint - \textit{contrast-preserving} signal and image restoration. In contrast-preserving restoration, the recovered signal or image is expected to keep not only neat edges, but also the contrasts of the edges. This is the perfect target in signal and image restoration problems. Some edge preserving models, such as TV models, still suffer from contrast reduction effect. This paper shows that any convex or smooth regularization (satisfying some common assumptions) based variational model is impossible or with low probabilities to preserve edge contrasts. This generalizes Meyer's remarkable example to any convex regularization and any clean data in discrete setting. Moreover, it is better to use those nonsmooth regularizers flat on $(\tau,+\infty)$ for some positive $\tau$, if to preserve the edge contrasts well. These discussions are new and further demonstrate the advantages of nonconvex and nonsmooth potential functions. They also naturally motivate us to present a general framework based on truncating existing potential functions. Some theoretical analysis in 1D and numerical algorithms with experiments in 2D for this framework show its better contrast preservation ability. Here we mention, truncated $\ell_p, 0<p\leq2$ exists in the literature with only numerical simulations. Also there is a very famous related regularizer called SCAD (smoothly clipped absolute deviation) \cite{Fan2001Variable}. It has closed form formula for the proximity operator, and has important applications in 1D
signal recovery~\cite{Fan2001Variable} and image deblurring with impulsive noise~\cite{Gu2017TVSCAD}. SCAD can be viewed as a truncation of $\ell_1$ if the middle quadratic part vanishes.

The paper is organized as follows. Section 2 gives the variational signal and image restoration models. We list some basic assumptions and the concept of sub-differential in section 3. In section 4, we discuss the motivation of our regularization framework. Section 5 is the general regularization framework based on truncation. In section 6, we give some theoretical analysis for 1D signal restoration using truncated regularization. In section 7, we present numerical methods to solve the 2D image restoration models with truncated regularization. Implementation details, convergence analysis and numerical experiments are included. Section 8 concludes the paper.

\section{The variational signal and image restoration models}

\subsection{1D signal}
We assume that a 1D signal $f\in \mathbb{R}^{K}$ is a degradation of $\underline{u}\in \mathbb{R}^N$:
$$
f=A\underline{u}+n,
$$
where $A$ is a matrix representing a linear operator such as a blur convolution; $n=\{n_i: 1\leq i\leq K\}$ is random noise. In this paper, we consider the most common and important case, where $\{n_i: 1\leq i\leq K\}$ are independently and identically distributed gaussian random variables, and more specifically, $n_i\thicksim \mathcal{N}(0,\sigma^2), 1\leq i\leq K$.

The restoration problem is to recover $\underline{u}$ from $f$. A widely used generic model is
\begin{equation}\label{eq-generic-model-1D}
\min_{u\in \mathbb{R}^N} \left\{
\begin{aligned}
E(u) &= \sum\limits_{1\leq i\leq N}\varphi((\nabla_x u)_i) + \frac{\alpha}{2}\|Au-f\|_{\mathbb{R}^K}^2\\
&= \sum\limits_{1\leq i\leq N}\rho(|(\nabla_x u)_i|) + \frac{\alpha}{2}\|Au-f\|_{\mathbb{R}^K}^2
\end{aligned}
\right\},
\end{equation}
with different choices of $\rho(\cdot)$. Here $\nabla_x$ is the forward difference operator with a specific boundary condition, e.g., the periodic or Neumann boundary condition.

\subsection{2D image}
Without loss of generality, we assume that a 2D image $\mathbf{f}\in \mathbb{R}^{K\times K}$ is a degradation of $\underline{\mathbf{u}}\in \mathbb{R}^{N\times N}$:
$$
\mathbf{f}=\mathbf{A}\underline{\mathbf{u}}+\mathbf{n},
$$
where $\mathbf{A}$ is a linear operator such as a blur convolution (note here we do not assume its matrix representation); $\mathbf{n}=\{\mathbf{n}_{i,j}, 1\leq i,j\leq K\}$ is random noise. Again here we consider Gaussian noise, where $\{\mathbf{n}_{i,j}\thicksim \mathcal{N}(0,\sigma^2), 1\leq i,j\leq K\}$ are independently and identically distributed.

In 2D there are two widely used generic minimization models to recover $\underline{\mathbf{u}}$ from $\mathbf{f}$. They are the following anisotropic model
\begin{equation}\label{eq-generic-model-2Dani}
\min_{\mathbf{u}\in \mathbb{R}^{N\times N}} \left\{
\begin{aligned}
E_{\mathrm{ani}}(\mathbf{u}) &= \sum\limits_{1\leq i,j\leq N}(\varphi((\nabla_x \mathbf{u})_{i,j})+\varphi((\nabla_y \mathbf{u})_{i,j})) + \frac{\alpha}{2}\|\mathbf{A}\mathbf{u}-\mathbf{f}\|_{\mathbb{R}^{K\times K}}^2\\
&= \sum\limits_{1\leq i,j\leq N}(\rho(|(\nabla_x \mathbf{u})_{i,j}|)+\rho(|(\nabla_y \mathbf{u})_{i,j}|)) + \frac{\alpha}{2}\|\mathbf{A}\mathbf{u}-\mathbf{f}\|_{\mathbb{R}^{K\times K}}^2\\
\end{aligned}
\right\},
\end{equation}
and the isotropic model
\begin{equation}\label{eq-generic-model-2Diso}
\min_{\mathbf{u}\in \mathbb{R}^{N\times N}} \left\{
\begin{aligned}
E_{\mathrm{iso}}(\mathbf{u}) &= \sum\limits_{1\leq i,j\leq N}\psi((\nabla_x \mathbf{u})_{i,j},(\nabla_y \mathbf{u})_{i,j}) + \frac{\alpha}{2}\|\mathbf{A}\mathbf{u}-\mathbf{f}\|_{\mathbb{R}^{K\times K}}^2\\
&= \sum\limits_{1\leq i,j\leq N}\rho(\sqrt{(\nabla_x \mathbf{u})_{i,j}^2+(\nabla_y \mathbf{u})_{i,j}^2}) + \frac{\alpha}{2}\|\mathbf{A}\mathbf{u}-\mathbf{f}\|_{\mathbb{R}^{K\times K}}^2\\
\end{aligned}
\right\},
\end{equation}
with different choices of $\rho(\cdot)$. Here $\nabla_x$ and $\nabla_y$ are the forward difference operators with a specific boundary condition. Sometimes we denote $\nabla=(\nabla_x,\nabla_y)$ for simplicity.


\section{Assumptions on $\rho(\cdot)$ and the sub-differential}
We make some assumptions about $\rho(\cdot)$ and give the notation of the sub-differential of a proper real function.

We first list some assumptions about $\rho(\cdot)$:
\begin{enumerate}
\item[(AS1)] $\rho(0)=0,\rho(s)<+\infty,\forall s<+\infty$ with $0$ as its strict minimizer;
\item[(AS2)] $\rho(s)$ is increasing over $[0,\infty)$;
\item[(AS3)] $\rho(s)$ is $C^2$ on $(0,+\infty)$;
\item[(AS4)] $\rho^{\prime\prime}(s) < 0$ strictly increases on $(0,+\infty)$ or $\rho^{\prime\prime}(s) \equiv 0$ on $(0,+\infty)$.
\end{enumerate}
Almost all regularizers in the image processing literature meet (AS1)(AS2).
Those regularizers non-differentiable at zero, except SCAD, satisfy (AS1)(AS2)(AS3)(AS4). SCAD satisfies (AS1)(AS2).


We now introduce the concept of limiting sub-differential (sub-differential for briefness) of a real proper function $h:\mathbb{R}^d\rightarrow\overline{\mathbb{R}}:=(-\infty,+\infty]$. Here that $h$ is proper means that it is finite somewhere and never equals $-\infty$.
The domain of $h$ is $\mathrm{dom}h=\{z\in\mathbb{R}^d:h(z)<+\infty\}$. The sub-differential of $h$ at $z\in\mathrm{dom}h$ is defined as
\begin{equation}\label{eq-subdifferential}
\partial h(z):=\left\{
\begin{aligned}
z_*\in\mathbb{R}^d: & \exists z^k\rightarrow z, z_*^k\rightarrow z_*, \mathrm{with} \\
&h(z^k)\rightarrow h(z), \liminf\limits_{v\rightarrow z^k}\frac{h(v)-h(z^k)-\langle z_*^k,v-z^k\rangle}{\|v-z^k\|}\geq0, \forall k\\
\end{aligned}
\right\}.
\end{equation}
For this definition, the following closedness property holds
\begin{equation}\label{eq-subdifferential-closedness}
\left\{z_*\in\mathbb{R}^d: \exists z^k\rightarrow z, z_*^k\rightarrow z_*, \mathrm{with~} h(z^k)\rightarrow h(z), z_*^k\in\partial h(z^k)\right\}\subset\partial h(z).
\end{equation}
We mention two special cases of the above limiting sub-differential. For a convex function, the limiting sub-differential in \cref{eq-subdifferential} reduces to the classical sub-differential in convex analysis. One useful property of the sub-differential of a convex function, which does not hold by the limiting sub-differential of a nonconvex function, is the monotonicity. For a continuously differentiable function, the limiting sub-differential is nothing but the gradient.

The sub-differential \cref{eq-subdifferential} can be used to characterize a extreme point of a function. If $z$ is a such point of $h$, then
\begin{equation}
\label{eq-firstOrderOptimCondi}
\partial h(z)\ni 0.
\end{equation}
The reader is referred to \cite{Rockafeller1998Variational} for more details.

Note $\varphi(x)=\rho(|x|)$ and $\psi(x,y)=\rho(\sqrt{x^2+y^2})$. The symmetry of $\varphi(\cdot)$ and $\psi(\cdot,\cdot)$ indicates the following basic properties on their sub-differentials.
\begin{proposition}
\label{proposition-subdiff}
Assume $\rho(\cdot)$ to satisfy (AS1)(AS2)(AS3). Then $\partial\varphi(x)$ and $\partial\psi(x,y)$ satisfy
\begin{enumerate}
\item[(1)] $g\in\partial\varphi(x)\Longleftrightarrow-g\in\partial\varphi(-x)$;
\item[(2)] $\partial\varphi(0)$ is either a closed interval $[-\rho\prime(0^+),\rho\prime(0^+)]$ or $(-\infty,+\infty)$;
\item[(3)] $(g^1,g^2)\in\partial\psi(x,y)\Longleftrightarrow (g^1,-g^2)\in\partial\psi(x,-y)$;\\
      $(g^1,g^2)\in\partial\psi(x,y)\Longleftrightarrow (-g^1,g^2)\in\partial\psi(-x,y)$;\\
      $(g^1,g^2)\in\partial\psi(x,y)\Longleftrightarrow (g^2,g^1)\in\partial\psi(y,x)$;
\item[(4)] $\partial\psi(0,0)$ is either a closed disk with radius $\rho\prime(0^+)$ centered at the origin or the whole plane $\mathbb{R}^2$;
\end{enumerate}
\end{proposition}
\proof
The symmetry of $\partial\varphi$ and $\partial\psi$ in (1)(3) is obvious. To show (2)(4), we notice that $\rho\prime(0^+)=\lim\limits_{t\rightarrow0^+}\frac{\rho(t)-\rho(0)}{t}=\lim\limits_{t\rightarrow0^+}\rho\prime(\theta(t) t)=\lim\limits_{s\rightarrow0^+}\rho\prime(s)$, where $0\leq\theta(t)\leq1,\forall t$. Then (2)(4) follow by the definition of the sub-differential in \cref{eq-subdifferential}.
\hfill$\blacksquare$

\section{Motivation}

It is well known that edges and contrasts are the most important features of signals and images.
An good method in signal and image restoration is capable to preserve both the edges and contrasts.
For variational methods, we first have two propositions on the perfect recovery for 1D signal, showing that nonconvex and nonsmooth regularization is desired.

The first one indicates that \cref{eq-generic-model-1D} with any convex regularizer is impossible to perfectly recover a nonconstant signal, even in the ideal case of $A=I$ and noiseless observed data.
\begin{proposition}\label{theorem-1Dgeneral-clean-convex}
Assume $\rho(\cdot)$ to be convex and satisfy (AS1)(AS2). If a signal $\tilde{u}\in\mathbb{R}^{N}$ can be recovered by the minimization problem \eqref{eq-generic-model-1D} with $f=A\tilde{u}$, then $\tilde{u}\in\mathbb{R}^{N}$ is a constant signal, i.e., $\tilde{u}=c(1,1,\cdots,1)\in\mathbb{R}^{N}$ for some $c\in \mathbb{R}$.
\end{proposition}
\proof
Denote $R(\tilde{u})=\sum\limits_{1\leq i\leq N}\varphi((\nabla_x \tilde{u})_i)$.
Since $\rho(\cdot)$ is convex, $R(\cdot)$ is also convex.
If $\tilde{u}\in\arg\min\limits_{u}E(u)$, we have
$$
\partial R(\tilde{u})+
\alpha A^{\textsf{T}}(A\tilde{u}-f)\ni0,
$$
where $\textsf{T}$ denotes the adjoint operation.
As $f=A\tilde{u}$ is noise free, the above becomes
$$
\partial R(\tilde{u})\ni0.
$$
Using the non-negativeness of $R(\tilde{u})$ and $R(0)-R(\tilde{u})\geq \langle0-\tilde{u},\partial R(\tilde{u})\rangle$, we get $R(\tilde{u})=0$. Since $0$ is the strict minimizer of the regularization function $\rho$, there exists $c\in \mathbb{R}$, such that $\tilde{u}=c(1,1,\cdots,1)\in\mathbb{R}^{N}$.
\hfill$\blacksquare$

When the observation $f$ is contaminated with random noise, we have an upper bound of the probability of the perfect recovery by \eqref{eq-generic-model-1D}.
\begin{proposition}\label{theorem-1Dgeneral-noise}
Assume $\rho$ to satisfy (AS1)(AS2) and $\nabla_x$ to be with Neumann boundary condition. Consider a general signal $\tilde{u}\in\mathbb{R}^N$ and its noisy observation $f=A\tilde{u}+n$.
For $i=1,\cdots,N$, denote
\begin{equation}
\label{eq-1Dgeneral-noise}
 \mathcal{P}_i = \Phi\left(\frac{\sup\{g: g\in \partial\varphi((\nabla_x\tilde{u})_{i-1})-\partial\varphi((\nabla_x\tilde{u})_{i})\}}{\alpha \sigma|a_i|_2}\right) - \Phi\left(\frac{\inf\{g: g\in \partial\varphi((\nabla_x\tilde{u})_{i-1})-\partial\varphi((\nabla_x\tilde{u})_{i})\}}{\alpha \sigma|a_i|_2}\right),
\end{equation}
with $\Phi(\cdot)$ as the cumulative distribution function of the standard normal distribution; and $a_i$ as the $i$th column of $A$.
Then we have
\begin{itemize}
\item[(1)] $\mathcal{P}(\tilde{u}\in \arg\min\limits_{u}E(u)) \leq \min\{\mathcal{P}_i,i=1,\cdots,N\}$;
\item[(2)] If $A=I$ (denoising problem), $\mathcal{P}(\tilde{u}\in \arg\min\limits_{u}E(u)) \leq \prod\limits_{i=1}^{N}\mathcal{P}_i$ with $|a_i|_2=1, \forall i$.
\end{itemize}
In addition, if the regularization function $\varphi(\cdot)$ is continuously differentiable everywhere, the recovery probability $\mathcal{P}(\tilde{u}\in \arg\min\limits_{u}E(u))=0$.
\end{proposition}
\proof The operator $\nabla_x$ with the Neumann boundary condition is surjective. Therefore the chain rule of sub-differential holds \cite{Rockafeller1998Variational}. Assume $\tilde{u}\in \arg\min\limits_{u}E(u)$. It then follows, by \cref{eq-firstOrderOptimCondi}, that
$$
\left(\nabla_x^{\textsf{T}}\prod\limits_{k=1}^{N}\partial\varphi((\nabla_x\tilde{u})_{k})\right)_{i}+
\alpha (A^{\textsf{T}}(A\tilde{u}-f))_i\ni0,\ 1\leq i\leq N,
$$
where $\partial\varphi((\nabla_x\tilde{u})_{k})$ is the sub-differential of $\varphi$ at $(\nabla_x\tilde{u})_{k}$; and $\textsf{T}$ denotes the adjoint operation. Denote the sub-gradients in $\partial\varphi((\nabla_x\tilde{u})_{k})$ as $g_{k}$, $k=1,\cdots,N$. Then,
\begin{equation*}
\begin{split}
&\phantom{\;=\;} \mathcal{P}(\tilde{u} \in \arg\min\limits_{u}E(u))\\
&= \mathcal{P}(\exists g_{i}\in\partial\varphi((\nabla_x\tilde{u})_{i}), 1\leq i\leq N, \mathrm{s.t.}, g_{i-1}-g_{i} = \alpha  (A^{\textsf{T}}n)_{i}, 1\leq i\leq N)\\
&\leq\mathcal{P}(\forall 1\leq i\leq N, \exists g_{i}\in\partial\varphi((\nabla_x\tilde{u})_{i}), g_{i-1}\in\partial\varphi((\nabla_x\tilde{u})_{i-1}), \mathrm{s.t.}, g_{i-1}-g_{i}= \alpha  (A^{\textsf{T}}n)_{i})	
\end{split}		
\end{equation*}
due to the randomness of $n$.

If we introduce, for each $i,1\leq i\leq N$, the random event $\mathcal{E}_{i}$ as
\begin{equation}\label{eq-1Dgeneral-randomevent1}
\partial\varphi((\nabla_x\tilde{u})_{i-1}) - \partial\varphi((\nabla_x\tilde{u})_{i}) \ni \alpha  (A^{\textsf{T}}n)_i,\ 1\leq i\leq N,
\end{equation}
it then follows that
\begin{equation}\label{eq-1Dgeneral-noise-prob}
\mathcal{P}(\tilde{u} \in \arg\min\limits_{u}E(u))\leq\mathcal{P}\left(\bigcap\limits_{i=1}^{N}\mathcal{E}_{i}\right).
\end{equation}
Since $n_j\thicksim \mathcal{N}(0,\sigma^2), 1\leq j\leq N$ are i.i.d, it follows that
$$
\mathcal{P}(\tilde{u}\in \arg\min\limits_{u}E(u)) \leq\min\{\mathcal{P}_i,i=1,\cdots,N\},
$$
which is exactly (1). Here $\mathcal{P}_i$ is defined in \cref{eq-1Dgeneral-noise}.
If $A=I$ (i.e., denoising problem), then the random events $\mathcal{E}_{i}$ defined in \eqref{eq-1Dgeneral-randomevent1} are independent. Therefore (2) holds. \hfill
$\blacksquare$

\begin{remark}\label{rem-1D-noise-prob}
For a piecewise constant signal, there are more samples in flat regions (with zero gradient) than discontinuities. To enlarge the upper bound of recovery probability, it is better to have $\rho\prime(0^+)$ as large as possible, since $\partial\varphi(0)-\partial\varphi(0)=[-2\rho\prime(0^+),2\rho\prime(0^+)]$.
\end{remark}

Although nonconvex and nonsmooth potential functions have attracted much attention due to its capability to produce high quality restorations with neat edges, the above two results further indicate the merits to use them.

Unfortunately, most nonconvex and nonsmooth potential functions still yield a contrast reduction effect, even in the ideal case of noiseless observation.
To see this clearly, let us consider the restoration of a typical signal by \cref{eq-generic-model-1D}. With $N=3M$, we define the following discrete gate function
\begin{equation}\label{eq-example-1D-clean}
\underline{u}_i=\left\{
\begin{array}{lll}
0,           \qquad 1 \leq i \leq M,\\
\bar{U},     \qquad M+1 \leq i \leq 2M,\\
0,           \qquad 2M+1\leq i \leq 3M,
\end{array}
\right.
\end{equation}
where $\bar{U}$ ($\bar{U}>0$, without loss of generality) is a constant. Obviously, $\underline{u}$ has two edge points and the contrasts are both $\bar{U}$.

\begin{proposition}\label{theorem-1D-clean-lowercontr}
Assume that $\rho(\cdot)$ satisfies (AS1)(AS2) and $\rho^{\prime}(\bar{U}^+)>0,\rho^{\prime}(\bar{U}^-)>0$ exist. Then
\begin{itemize}
\item[(1)] Any $u$ with higher contrast than $\underline{u}$ satisfies $E(u)>E(\underline{u})$;
\item[(2)] There always exists $u$ with lower contrast than $\underline{u}$ satisfying $E(u)<E(\underline{u})$.
\end{itemize}
\end{proposition}
\proof We prove the result by contradiction. \\
(1) If there exists a $\hat{u}\in \mathbb{R}^N$ with higher contrast than $\underline{u}$ satisfying $E(\hat{u})\leq E(\underline{u})$, then $|\hat{u}_{M+1}-\hat{u}_{M}|\geq\bar{U}, |\hat{u}_{2M+1}-\hat{u}_{2M}|\geq\bar{U}$ with at least one strict inequality, and the following holds
$$
0\geq E(\hat{u})-E(\underline{u})\geq \rho(|\hat{u}_{M+1}-\hat{u}_{M}|)+\rho(|\hat{u}_{2M+1}-\hat{u}_{2M}|)-2\rho(\bar{U})\geq0,
$$
due to the monotonicity of $\rho$. It is necessary that $\rho(|\hat{u}_{M+1}-\hat{u}_{M}|)=\rho(|\hat{u}_{2M+1}-\hat{u}_{2M}|)=\rho(\bar{U})$. Since either $|\hat{u}_{M+1}-\hat{u}_{M}|>\bar{U}$ or $|\hat{u}_{2M+1}-\hat{u}_{2M}|>\bar{U}$, it follows that $\rho^{\prime}(\bar{U}^+)=0$. This is a contradiction.\\
(2) Let us define $u^{\epsilon}\in \mathbb{R}^N$ as follows
$$
u^{\epsilon}_i=\left\{
\begin{array}{lll}
\epsilon,           \qquad 1 \leq i \leq M,\\
\bar{U}-\epsilon,     \qquad M+1 \leq i \leq 2M,\\
\epsilon,           \qquad 2M+1\leq i \leq 3M,
\end{array}
\right.
$$
where $\epsilon>0$ is a small positive number. Obviously $u^{\epsilon}$ has lower contrast than $\underline{u}$.
Now we show that there exists an $\epsilon>0$ satisfying $E(u^{\epsilon})<E(\underline{u})$.
Assume that $E(u^{\epsilon})\geq E(\underline{u}), \forall \epsilon>0$. That is,
$$
2\rho\left(\bar{U}-2\epsilon\right) - 2\rho\left(\bar{U}\right) + \frac{\alpha}{2}\|A(\underline{u}-\epsilon\mathbf{1}_{\mathrm{mid}}+\epsilon(\mathbf{1}-\mathbf{1}_{\mathrm{mid}}))-f\|_{\mathbb{R}^{K}}^2\geq0,\ \forall \epsilon>0,
$$
where $\mathbf{1}=(1,1,\cdots,1)\in\mathbb{R}^{N}$, $\mathbf{1}_{\mathrm{mid}}=(\underbrace{0,\cdots,0}\limits_M,\underbrace{1,\cdots,1}\limits_M,\underbrace{0,\cdots,0}\limits_M)\in\mathbb{R}^{N}$. Then
$$
2\rho\left(\bar{U}-2\epsilon\right) - 2\rho\left(\bar{U}\right) + \epsilon^2\frac{\alpha}{2}\|A(-\mathbf{1}_{\mathrm{mid}}+(\mathbf{1}-\mathbf{1}_{\mathrm{mid}}))\|_{\mathbb{R}^{K}}^2\geq0,\ \forall \epsilon>0.
$$
Letting $\epsilon\rightarrow0$, we get
$$
\rho^{\prime}\left(\bar{U}^-\right)\leq0,
$$
which is a contradiction.
\hfill$\blacksquare$

\begin{remark}\hfill
\begin{itemize}
\item Most nonconvex potential functions in the literature satisfy the conditions of Proposition~\ref{theorem-1D-clean-lowercontr}.
\item The above proposition can be extended to more general signals than \cref{eq-example-1D-clean}.
\end{itemize}
\end{remark}

\begin{remark}
By the monotonicity of the sub-differential of convex functions, $\rho^{\prime}(\bar{U}^+)>0$ and $\rho^{\prime}(\bar{U}^-)>0$ hold for any convex $\rho(\cdot)$ with (AS1). Therefore, \cref{eq-generic-model-1D} with any convex $\rho(\cdot)$ satisfying (AS1) has contrast reduction effect.
\end{remark}

These discussions can be extended to 2D image case. We omit the details. As can be seen, it is reasonable to consider those potential functions flat on $(\tau,+\infty)$ for some positive $\tau$ related to the edge contrasts in variational signal and image restoration. This motivates us to truncate existing potential functions (especially those non-differentiable at $0$), to generate and use new potential functions.

\section{Truncated regularization: a general regularization framework}

For a general regularizer function $\rho(\cdot)$, we construct a new regularizer function
\begin{equation}
\label{eq-TR-regularizer}
\T (\cdot)=\rho_{\tau}(\cdot)=\rho(\min(\cdot,\tau)),
\end{equation}
where $\tau>0$ is a positive real parameter. This regularizer function is flat on $(\tau,+\infty)$. It is easy to see that, if $\rho(\cdot)$ satisfies the basic assumptions (AS1)(AS2), $\T(\cdot)=\rho_{\tau}(\cdot)$ also satisfies the basic assumptions (AS1)(AS2). However, $\T(\cdot)$ is always nonconvex. For $\rho_{5}(\cdot)$ in \cref{tab-PFs}, its truncated version is the same as itself. For a general regularizer $\rho(\cdot)$, $\T (\cdot)=\rho_{\tau}(\cdot)$ degenerates to $\rho(\cdot)$, when $\tau=+\infty$. In the following, $\tau$ takes a finite positive real value, unless pointed out specifically. Note that $\rho_{8}(\cdot)$ (SCAD) in \cref{tab-PFs} is already a smoothed and more general version of truncated $\ell_1$; see Remark \ref{rem-SCAD}.


Basic calculations show that, for those regularizer functions in the literature with non-differentiability at $0$ (e.g. the regularizer functions in \cref{tab-PFs}), the so called subadditivity holds. Actually, their truncated versions (SCAD is already a smoothed truncation of $\ell_1$) also have this property.
\begin{lemma}\label{subadditivity-minfunction}
Given $a,b\geq 0, \tau>0$, then
\begin{equation}
\min(a+b,\tau)\leq \min(a,\tau)+\min(b,\tau).
\end{equation}
\end{lemma}
\proof
Case 1: $0\leq a+b\leq \tau$. Then $0\leq a,b \leq \tau$. Therefore, $\min(a+b,\tau)=a+b=\min(a,\tau)+\min(b,\tau)$.

Case 2: $a+b>\tau$. If $0\leq a,b<\tau$, then $\min(a,\tau)+\min(b,\tau)=a+b>\tau=\min(a+b,\tau)$. If $a\geq \tau$ or $b \geq \tau$, then
$\min(a,\tau)+\min(b,\tau)\geq \tau=\min(a+b,\tau)$.\hfill$\blacksquare$

\begin{proposition}\label{subadditivity-TRfunction}
Given $\tau>0$, if $\rho(\cdot)$ satisfies the subadditivity property over $[0,+\infty)$ and the assumptions (AS1)(AS2), then its truncated version $\T(\cdot)=\rho(\min(\cdot,\tau))$ also has the subadditivity property over $[0,+\infty)$.
\end{proposition}
\proof For $\forall a,b\geq 0$, $\T(a+b)=\rho(\min(a+b,\tau))\leq\rho(\min(a,\tau)+\min(b,\tau))
\leq\rho(\min(a,\tau))+\rho(\min(b,\tau))=\T(a)+\T(b)$.
\hfill $\blacksquare$

\begin{remark}\label{rem-SCAD}
SCAD uses a sophisticated quadratic function to smooth the truncated $\ell_1$ at the truncation point locally. When the middle quadratic part vanishes, SCAD degenerates to truncated $\ell_1$. Therefore, it is also a more general version of truncated $\ell_1$. By its concavity on $[0,+\infty)$, it is easy to check that SCAD also has the subadditivity property.
\end{remark}

\begin{table}[htbp]
  \caption{Some potential functions satisfying the subadditivity, where $0<p<1, \theta > 0$,  $a > 2$}
  \label{tab-PFs}
  \centering
  \begin{tabular}{p{0.1cm} l p{0.1cm}} \toprule
    & $\rho_1(s) = s$                                       & \\
    & $\rho_2(s) = s^p$                                     & \\
    & $\rho_3(s) = \ln(\theta s + 1) $                    & \\
    & $\rho_4(s) = \frac{\theta s}{ 1 + \theta s}$      & \\
    & $\rho_5(0) = 0, \rho_5(s) = 1$ if $s > 0 $      & \\
    & $\rho_6(s) = \ln(\theta s^p + 1)$                   & \\
    & $\rho_7(s) = \frac{\theta s^p}{ 1 + \theta s^p }$ & \\
    & $\rho_8(s)=
    \left\{
      \begin{array}{ll}
        \theta s, & s \leq \theta \\
        \frac{-s^{2} - \theta^{2} + 2 a \theta s}{2(a-1)}, & \theta < s < a \theta  \\
        \frac{(a+1)\theta^{2}}{2}, & s > a \theta
      \end{array}
    \right.$
    & \\
    \bottomrule
  \end{tabular}
\end{table}

With $\T(\cdot)=\rho(\min(\cdot,\tau))$ or the SCAD regularizer $\rho_8(\cdot)$ replacing $\rho(\cdot)$ in $\varphi$ and $\psi$ in the variational problems \eqref{eq-generic-model-1D} \eqref{eq-generic-model-2Dani} and \eqref{eq-generic-model-2Diso}, the energy functions $E(\cdot)$ in \eqref{eq-generic-model-1D}, $E_{\mathrm{ani}}(\cdot)$ in \eqref{eq-generic-model-2Dani}, and $E_{\mathrm{iso}}(\cdot)$ in \eqref{eq-generic-model-2Diso} are coercive if and only if $A^{\textsf{T}}A$ ($\mathbf{A}^{\textsf{T}}\mathbf{A}$) is invertible. Therefore, we assume:
\begin{itemize}
\item[(AS5)] $A^{\textsf{T}}A$ ($\mathbf{A}^{\textsf{T}}\mathbf{A}$) is invertible.
\end{itemize}

\section{Truncated regularization for 1D signal restoration: Some analytic results}
In this section we prove some analytic results on variational 1D signal restoration using truncated regularization, showing its better contrast-preserving ability.

Assume $\emptyset\neq\Omega\subsetneq J = \{1,\cdots,N\}$. Let $\mathbf{1}_{\Omega}$ be its indicator function and $\zeta>0$ be a real number. We define two index sets
$$
J_0=\{i:(\nabla_x\mathbf{1}_{\Omega})_i=0\};\ J_1=\{i:(\nabla_x\mathbf{1}_{\Omega})_i\neq0\}= J\setminus J_0.
$$
Consider now the minimization problem \eqref{eq-generic-model-1D} using truncated regularization where $f=A(\zeta\mathbf{1}_{\Omega})\in \mathbb{R}^{K}$.
For clarity, let us denote
\begin{equation}\label{eq-1D-TR-model}
E^{\zeta}(u) = \sum\limits_{1\leq i\leq N}\T(|(\nabla_x u)_i|) + \frac{\alpha}{2}\|A(u-\zeta\mathbf{1}_{\Omega})\|_{\mathbb{R}^K}^2,
\end{equation}
where the regularizer function $\T$ reads \eqref{eq-TR-regularizer}.

The following theorem shows the perfect recovery (i.e., contrast preservation) of the signal $\zeta\mathbf{1}_{\Omega}$ by \cref{eq-generic-model-1D} with a truncated regularization. The proof is motivated by \cite{Nikolova2005Analysis}.
\begin{theorem}\label{theorem-1D-TR-globalminimizer}
If $\zeta>\tau+\sqrt{\frac{4\T(\tau)}{\alpha\mu_{\min}}\#J_1}$, then the global minimizer is $\zeta\mathbf{1}_{\Omega}$. Here $\mu_{\min}>0$ is the minimal eigenvalue of $A^{\textsf{T}}A$.
\end{theorem}
\proof
It is clear that $E^\zeta(\zeta\mathbf{1}_{\Omega})=\sum\limits_{1\leq i\leq N}\T(|(\nabla_x \zeta\mathbf{1}_{\Omega})_i|)=\sum\limits_{i\in J_1}\T(\zeta)\leq\#J_1\T(\tau)$.

Assume $v$ to be the global minimizer. We define two index sets
$$
J^v_{0}=\{i:|(\nabla_xv)_i|\leq \tau\};\ J^v_{1}=\{i:|(\nabla_xv)_i|> \tau\}= J\setminus J^v_0.
$$
We first check $J^v_0=J_0;J^v_1=J_1$ by contradiction. Assume $J^v_0\cap J_1\neq\emptyset$. There exists $i$ such that
$$
|(\nabla_xv)_i|=|v_{i+1}-v_{i}|\leq \tau;\ |(\nabla_x\zeta\mathbf{1}_{\Omega})_i|=\zeta.
$$
It then follows, under the assumption $\zeta>\tau+\sqrt{\frac{4\T(\tau)}{\alpha\mu_{\min}}\#J_1}$, that
\begin{align*}
E^{\zeta}(v)\geq\frac{\alpha}{2}\|A(v-\zeta\mathbf{1}_{\Omega})\|_{\mathbb{R}^K}^2&\geq\frac{\alpha}{2}\mu_{\min}\|v-\zeta\mathbf{1}_{\Omega}\|_{\mathbb{R}^K}^2\\
&\geq\frac{\alpha}{2}\mu_{\min}((v_{i+1}-\zeta)^2+(v_{i}-0)^2)\\
&\geq\frac{\alpha}{2}\mu_{\min}\frac{1}{2}(\zeta-\tau)^2\\
&>\#J_1\T(\tau)\geq E^{\zeta}(\zeta\mathbf{1}_{\Omega}),
\end{align*}
which is a contradiction. It holds that $J^v_0\cap J_1=\emptyset$. Now assume $J_1\subsetneq J^v_1$. Then $\#J^v_1\geq\#J_1+1$ and $\|A(v-\zeta\mathbf{1}_{\Omega})\|_{\mathbb{R}^K}^2>0$. By the definition of $J^v_1$,
it then follows that
$$
E^{\zeta}(v)>\#J^v_1\T(\tau)\geq(\#J_1+1)\T(\tau)>(\#J_1)\T(\tau)\geq E^{\zeta}(\zeta\mathbf{1}_{\Omega}),
$$
which is a contradiction. Therefore, $J^v_0=J_0;J^v_1=J_1$.

If the global minimizer $v\neq \zeta\mathbf{1}_{\Omega}$, then $\|A(v-\zeta\mathbf{1}_{\Omega})\|_{\mathbb{R}^K}^2>0$ and
$$
\sum\limits_{i\in J^v_1}\T(|(\nabla_xv)_i|)=\sum\limits_{i\in J^v_1}\T(\tau)=\sum\limits_{i\in J_1}\T(\zeta)=\sum\limits_{i\in J_1}\T(|(\nabla_x\zeta\mathbf{1}_{\Omega})_i|),
$$
where we used $\zeta > \tau$ by the assumption. It then follows that
$$
E^{\zeta}(v)>E^{\zeta}(\zeta\mathbf{1}_{\Omega}),
$$
which is a contradiction. Therefore, the global minimizer is $\zeta\mathbf{1}_{\Omega}$.
\hfill$\blacksquare$

\begin{remark}\label{rem-1D-TR-perfectRecovery}\hfill
\begin{itemize}
\item From the proof, one can see that \cref{theorem-1D-TR-globalminimizer} holds for any truncated regularizers.
\item \cref{theorem-1D-TR-globalminimizer} can be extended to 2D case.
\end{itemize}
\end{remark}

For a general parameter setting, it is quite difficult to give the structure of the global minimizer. However, if the regularizer function satisfies the subadditivity and the $A^{\textsf{T}}A$ is a diagonal matrix (with image denoising as a special case where $A^{\textsf{T}}A=I$), we can show some structure of the global minimizer.

\begin{theorem} \label{theorem-1D-TR-jumpingpointset}
Assume $A^{\textsf{T}}A = \diag\{d_1,d_2,\cdots,d_N\}, d_i > 0, i = 1,2,\cdots,N$ and the regularizer function $\T(\cdot)$ in \eqref{eq-1D-TR-model} to satisfy the subadditivity (Note here we do not require the finiteness of $\tau$). Let $v$ be the global minimizer of \eqref{eq-1D-TR-model}, then
\begin{itemize}
\item[(1)] The extremum principle holds, i.e., $0\leq v_i\leq \zeta$ for all $i \in J$;
\item[(2)] No new (and thus false) discontinuity appears in $v$, i.e., $J^*_{1}=\{i\in J:(\nabla_x v)_i\neq0\} \subseteq J_1$;
\item[(3)] $v$ preserves the monotonicity (not necessarily strict) from the input signal $\zeta\mathbf{1}_{\Omega}$.
\end{itemize}
\end{theorem}

\proof (1) We show, by contradiction, that $0\leq v_i\leq \zeta$ for all $i \in J$. Assume that $v_j < 0$ for some $j \in J$. We define a new signal $\bar{v}$ by
\begin{displaymath}
  \bar{v}_i=\left\{
  \begin{array}{ll}
    0,   & i \in \{j:v_j < 0\} ,\\
    v_i, & \mathrm{otherwise.}
  \end{array}
  \right.
\end{displaymath}
For each $i$, there are four cases for $|v_{i+1}-v_{i}|$ and $|\bar{v}_{i+1}-\bar{v}_{i}|$:
\begin{enumerate}
\item ${v}_{i+1},{v}_{i} \geq 0$:
      $|\bar{v}_{i+1}-\bar{v}_{i}|=|{v}_{i+1}-{v}_{i}|$;
\item ${v}_{i+1},{v}_{i} < 0$:
      $|\bar{v}_{i+1}-\bar{v}_{i}|=0\leq |{v}_{i+1}-{v}_{i}|$;
\item ${v}_{i+1} \geq 0, {v}_{i} <0$:
      $|\bar{v}_{i+1}-\bar{v}_{i}|=|v_{i+1}-0|=v_{i+1} < v_{i+1}-v_i= |{v}_{i+1}-{v}_{i}|$;
\item ${v}_{i+1} < 0, {v}_{i} \geq 0$:
      $|\bar{v}_{i+1}-\bar{v}_{i}|=|0-v_{i}|=v_{i} < v_{i}-v_{i+1}= |{v}_{i+1}-{v}_{i}|$.
\end{enumerate}
It follows that $|(\nabla_x v)_{i}| = |v_{i+1}-v_{i}| \geq |\bar{v}_{i+1}-\bar{v}_{i}| = |(\nabla_x \bar{v})_{i}|$ always holds. Thus,
$$
  \T(|(\nabla_x v)_{i}|) \geq \T(|(\nabla_x \bar{v})_{i}|), \; i \in J.
$$
It implies that
\begin{equation*}
  \begin{split}
     E^{\zeta}(v) - E^{\zeta}(\bar{v})
     &= \sum\limits_{i \in J}\T(|(\nabla_x v)_i|)
      - \sum\limits_{i \in J}\T(|(\nabla_x \bar{v})_i|)
      + \frac{\alpha}{2}\|A(v-\zeta\mathbf{1}_{\Omega})\|_{\mathbb{R}^K}^2-
      \frac{\alpha}{2}\|A(\bar{v}-\zeta\mathbf{1}_{\Omega})\|_{\mathbb{R}^K}^2\\
     &\geq       \frac{\alpha}{2}\|A(v-\zeta\mathbf{1}_{\Omega})\|_{\mathbb{R}^K}^2-
      \frac{\alpha}{2}\|A(\bar{v}-\zeta\mathbf{1}_{\Omega})\|_{\mathbb{R}^K}^2\\
     &= \frac{\alpha}{2}\sum\limits_{i \in J} d_i |(v-\zeta\mathbf{1}_{\Omega})_i|^2
     - \frac{\alpha}{2}\sum\limits_{i \in J} d_i |(\bar{v}-\zeta\mathbf{1}_{\Omega})_i|^2\\
     &= \frac{\alpha}{2} \sum\limits_{i,v_i<0} d_i(|(v-\zeta\mathbf{1}_{\Omega})_i|^2
     -|(0-\zeta\mathbf{1}_{\Omega})_i|^2)>0.
  \end{split}
\end{equation*}
This contradicts the fact that $v$ is the global minimizer of \eqref{eq-1D-TR-model}. Thus we have $v_i\geq 0$ for all $i \in J$. In a similar way we can prove that $v_i \leq \zeta$ for all $i \in J$.

(2) Now we prove $J^*_{1}=\{i\in J:(\nabla_x v)_i\neq0\} \subseteq J_1$. If $J^*_1 = \emptyset$, then $J^*_1 \subseteq J_1$ holds. We consider that $J^*_1$ is nonempty. Suppose that there exists $j \in J^*_1$ but $j \notin J_1$. That is, $v_j \neq v_{j+1}$ and $j \in J_0$. In total there are four cases:
\begin{enumerate}
\item $v_j > v_{j+1}$ and $(\zeta\mathbf{1}_{\Omega})_j=(\zeta\mathbf{1}_{\Omega})_{j+1}=0$;
\item $v_j > v_{j+1}$ and $(\zeta\mathbf{1}_{\Omega})_j=(\zeta\mathbf{1}_{\Omega})_{j+1}=\zeta$;
\item $v_j < v_{j+1}$ and $(\zeta\mathbf{1}_{\Omega})_j=(\zeta\mathbf{1}_{\Omega})_{j+1}=0$;
\item $v_j < v_{j+1}$ and $(\zeta\mathbf{1}_{\Omega})_j=(\zeta\mathbf{1}_{\Omega})_{j+1}=\zeta$.
\end{enumerate}
Without loss of generality, we consider case (i). Define a new vector by
\begin{displaymath}
  \tilde{v}_i=\left\{
  \begin{array}{ll}
    v_{i + 1}, & i = j, \\
    v_i, & \mathrm{otherwise}.
  \end{array}
  \right.
\end{displaymath}
where $\tilde{v}$ only differs from $v$ at index $j$. We now show $E^{\zeta}(\tilde{v}) < E^{\zeta}(v)$.  Indeed, we have
{\small
\begin{equation*}
  \begin{split}
    &\phantom{\;=\;} E^{\zeta}(v)-E^{\zeta}(\tilde{v}) \\
    & = \sum\limits_{i \in J}\T(|(\nabla_x v)_i|)
      - \sum\limits_{i \in J}\T(|(\nabla_x \tilde{v})_i|)
      + \frac{\alpha}{2}\|A(v-\zeta\mathbf{1}_{\Omega})\|_{\mathbb{R}^K}^2-
      \frac{\alpha}{2}\|A(\tilde{v}-\zeta\mathbf{1}_{\Omega})\|_{\mathbb{R}^K}^2\\
    & = \T(|(\nabla_x v)_{j -1}|) + \T(|(\nabla_x v)_j|)
         -\T(|(\nabla_x \tilde{v})_{j-1}|)-\T(|(\nabla_x \tilde{v})_j|)
         + \frac{\alpha}{2} d_j|(v-\zeta\mathbf{1}_{\Omega})_j|^2
         - \frac{\alpha}{2} d_j|(\tilde{v}-\zeta\mathbf{1}_{\Omega})_j|^2\\
    & =  \T(|(\nabla_x v)_{j -1}|) + \T(|(\nabla_x v)_j|)
        -\T(|(\nabla_x \tilde{v})_{j-1}|)
      +  \frac{\alpha}{2} d_j(|v_j|^2 - |v_{j+1}|^2)\\
    & >  \T(|v_{j}-v_{j-1}|) + \T(|v_{j+1}-v_{j}|)-\T(|v_{j+1}-v_{j-1}|).
  \end{split}
\end{equation*}
}%
The last inequality follows from $v_j > v_{j+1}$ and the extremum principle proved in (1). Since it follows that $\T(|v_{j+1}-v_{j-1}|)\leq \T(|v_{j}-v_{j-1}|) + \T(|v_{j+1}-v_{j}|)$ by the subadditivity, we have $E^{\zeta}(v)-E^{\zeta}(\tilde{v})>0$. This contradicts that $v$ is a global minimizer. For other three cases we can obtain similarly (by defining different new vectors) this contradiction. Therefore, $i \in J^*_1$ implies $i \in J_1$. That is, $J^*_1 \subseteq J_1$.

(3) We prove the result by contradiction. Without loss of generality, we assume there exists $j\in J$ such that $(\zeta\mathbf{1}_{\Omega})_j\leq(\zeta\mathbf{1}_{\Omega})_{j+1}$ but $v_j > v_{j+1}$. If $(\zeta\mathbf{1}_{\Omega})_j=(\zeta\mathbf{1}_{\Omega})_{j+1}$, (2) gives $v_j = v_{j+1}$, which is a contradiction. If $(\zeta\mathbf{1}_{\Omega})_j<(\zeta\mathbf{1}_{\Omega})_{j+1}$, by defining
\begin{displaymath}
  \hat{v}_i=\left\{
  \begin{array}{ll}
    v_{i + 1}, & i = j, \\
    v_i, & \mathrm{otherwise},
  \end{array}
  \right.
\end{displaymath}
and a similar argument in (2), we can show $E^{\zeta}(\hat{v}) < E^{\zeta}(v)$, which is a contradiction. Therefore, $v$ preserves the monotonicity (not necessarily strict by (2)) from the input $\zeta\mathbf{1}_{\Omega}$.
\hfill$\blacksquare$

\begin{remark}\label{rem-1D-SCAD-perfectRecovery}
From the proofs and the subadditivity property of SCAD, one can see that both \cref{theorem-1D-TR-globalminimizer} and \cref{theorem-1D-TR-jumpingpointset} hold for \eqref{eq-1D-TR-model} with $\T(\cdot)$ replaced by the SCAD regularizer.
\end{remark}


\section{Truncated regularization for 2D image restoration: Numerical method}
\label{sec-2D-TR}

We now present a numerical method with implementation details and convergence verification for variational 2D image restoration using truncated regularization.

For the simplicity of presentation, we denote the Euclidean space $\mathbb{R}^{K \times K}$ by $V$ and $(\mathbb{R}^{N \times N})\times (\mathbb{R}^{N \times N})$ by $Q$. The spaces $V$ and $Q$ are equipped with inner products $(\cdot, \cdot)_V,(\cdot, \cdot)_Q$ and norms $\|\cdot\|_V,\|\cdot\|_Q$, respectively; See, e.g., \cite{Wu2010Augmented}. In addition, we mention that for $\vecq\in Q$, $\vecq_{i,j}=(q^{(1)}_{i,j},q^{(2)}_{i,j})$ and
$|\vecq_{i,j}|=\sqrt{ (q^{(1)}_{i,j})^2 +( q^{(2)}_{i,j})^2 }$.

\subsection{Anisotropic and isotropic 2D truncated regularization models}

Our anisotropic and isotropic 2D truncated regularization models for image restoration are presented as follows:
\begin{itemize}
\item\textit{Anisotropic 2D truncated regularization model}:
\end{itemize}
\begin{equation}
\label{eq-2Dani-TR}
  \min\limits_{\mathbf{u}\in V}
  \left\{E_{\mathrm{ani}}^{\T}(\mathbf{u})=
  \sum_{1 \leq i,j \leq N}
  \T(|(\nabla_x \mathbf{u})_{ij}|) + \T(|(\nabla_y \mathbf{u})_{ij}|)
   + \frac{\alpha}{2} \|\mathbf{A} \mathbf{u} - \mathbf{f}\|_V^2\right\},
\end{equation}
\begin{itemize}
\item\textit{Isotropic 2D truncated regularization model}:
\end{itemize}
\begin{equation}
  \label{eq-2Diso-TR}
  \min\limits_{\mathbf{u}\in V}
  \left\{E_{\mathrm{iso}}^{\T}(\mathbf{u})=
  \sum_{1 \leq i,j \leq N}
  \T\left(\sqrt{(\nabla_x \mathbf{u})_{ij}^2 + (\nabla_{ y } \mathbf{u})_{ij}^2}\right)
  + \frac{\alpha}{2} \|\mathbf{A} \mathbf{u} - \mathbf{f}\|_V^2\right\},
\end{equation}
which are truncated versions of \eqref{eq-generic-model-2Dani} and \eqref{eq-generic-model-2Diso}, respectively.


\subsection{ADMM for 2D truncated regularization models}
\label{sec-admm}

We use the alternating direction method of multipliers (ADMM) to solve the truncated regularization models. ADMM is a very successful method to solve large-scale optimization problems with special structures. Since the methods to solve models \eqref{eq-2Dani-TR} and \eqref{eq-2Diso-TR} are very similar, we only elaborate the details for \eqref{eq-2Diso-TR} with the periodic boundary condition.


We first rewrite \eqref{eq-2Diso-TR} to the following constrained optimization problem:
\begin{equation}
\label{eq-2Diso-TR-equiv}
\begin{split}
    &   \min\limits_{(\mathbf{u},\vecq)\in V\times Q}
  \left\{{\tilde{E}}_{\mathrm{iso}}^{\T}(\mathbf{u},\vecq)=
      \sum\limits_{1 \leq i,j \leq N} \T(|\vecq_{ij}|)
      + \frac{\alpha}{2} \|\mathbf{A} \mathbf{u} - \mathbf{f}\|_V^2 \right\}, \\
    & \mathrm{s.t.}\quad\quad \vecq= (\nabla_x \mathbf{u}, \nabla_y \mathbf{u})
\end{split}
\end{equation}
and define the augmented Lagrangian functional for the problem \eqref{eq-2Diso-TR-equiv} as follows:
\begin{equation}
\label{eq-2Diso-TR-equiv-augLagfun}
\begin{split}
\mathcal{L}(\mathbf{u},\vecq;\vecmu)&=
\sum\limits_{1 \leq i,j \leq N} \T(|\vecq_{ij}|)
+\frac{\alpha}{2} \|\mathbf{A} \mathbf{u} - \mathbf{f}\|_V^2
+(\vecmu,\vecq-\nabla \mathbf{u})_Q +\frac{\beta}{2}\|\vecq-\nabla \mathbf{u}\|_Q^2\\
&=R(\vecq)
+\frac{\alpha}{2} \|\mathbf{A} \mathbf{u} - \mathbf{f}\|_V^2
+(\vecmu,\vecq-\nabla \mathbf{u})_Q +\frac{\beta}{2}\|\vecq-\nabla \mathbf{u}\|_Q^2,
\end{split}
\end{equation}
where $\nabla \mathbf{u} = (\nabla_x \mathbf{u}, \nabla_y \mathbf{u})$; $\vecmu\in Q$ is the Lagrangian multiplier; $\beta>0$ is a constant; and $R(\vecq)$ is introduced to simply the notation of the regularization term. The ADMM for solving \eqref{eq-2Diso-TR} can be described in \cref{alg-admm}.

\begin{algorithm}[t]
\caption{ADMM for the isotropic truncated regularization model \eqref{eq-2Diso-TR}}
\label{alg-admm}
\begin{algorithmic}[1]
  \STATE{Initialization: $\mathbf{u}^0,\vecq^0, \vecmu^0$;}
  \WHILE{stopping criteria is not satisfied}
  \STATE{Compute $\vecq^{k+1}$, $\mathbf{u}^{k+1}$, and update $\vecmu^{k+1}$ as follows:
  \begin{flalign}
   & \vecq^{k+1} \in \argmin_{\vecq \in Q} \mathcal{L}(\mathbf{u}^{k},\vecq;\vecmu^k), && \label{eq-admm-qsub}\\
   & \mathbf{u}^{k+1} = \argmin_{\mathbf{u}\in V} \mathcal{L}(\mathbf{u},\vecq^{k+1};\vecmu^k), && \label{eq-admm-usub}\\
   & \vecmu^{k+1}= \vecmu^{k} + \beta(\vecq^{k+1} - \nabla \mathbf{u}^{k+1}), && \label{eq-admm-lsub}
  \end{flalign}
  }
  \ENDWHILE
\end{algorithmic}
\end{algorithm}

In each iteration of ADMM, we need to solve the $\vecq-$sub problem \eqref{eq-admm-qsub} and the $\mathbf{u}-$sub problem \eqref{eq-admm-usub}.
The $\mathbf{u}-$sub problem \eqref{eq-admm-usub} is a quadratic optimization problem,
whose optimality condition gives a linear system
\begin{equation*}
\alpha \mathbf{A}^{\textsf{T}} (\mathbf{A} \mathbf{u}- \mathbf{f})-\nabla^{\textsf{T}}\vecmu^k - \beta\nabla^{\textsf{T}}\vecq^{k+1}- \beta\Delta \mathbf{u} =0,
\end{equation*}
which can be solved by Fourier transform with an FFT implementation \cite{Wang2008New,Wu2010Augmented}.

We now design a fast algorithm to find the global minimizer of the non-convex and nonsmooth $\vecq-$sub problem \eqref{eq-admm-qsub}.

\subsection{Find the global minimizer of the $\vecq-$sub problem \eqref{eq-admm-qsub}}
The problem \eqref{eq-admm-qsub} reads
$$
\min\limits_{\vecq\in Q} \left\{
     \sum\limits_{1 \leq i,j \leq N} \T(|\vecq_{ij}|)
     + (\vecmu^k,\vecq-\nabla \mathbf{u}^{k})_Q
     +\frac{\beta}{2}\|\vecq -\nabla \mathbf{u}^{k}\|_Q^2
     \right\},
$$
which, by the monotonicity of $\rho$ over $[0,+\infty)$, is
$$
     \min\limits_{\vecq\in Q} \left\{
     \sum\limits_{1 \leq i,j \leq N}
     \min\big(\rho(|\vecq_{ij}|), \rho(\tau)\big)
     +\frac{\beta}{2}|\vecq_{i,j} - \mathbf{w}_{i,j}|^2
     \right\},
$$
where $\mathbf{w}=\nabla \mathbf{u}^{k} - \vecmu^k /\beta \in Q$. This problem is separable. Consequently we only need to consider
\begin{equation*}
\min\limits_{\vecz\in \mathbb{R}^2}\left\{
g(\vecz;\vecw)=\min\big(\rho(|\vecz|), \rho(\tau)\big)
+\frac{\beta}{2}|\vecz - \vecw|^2\right\},
\end{equation*}
where $|\vecz|=\sqrt{(z^{(1)})^2 + (z^{(2)})^2}$; and $\vecw\in \R^2, \tau>0, \beta>0$ are given.

Suppose $\vecz^*=\argmin\limits_{\vecz \in \R^2}g(\vecz;\vecw)$. If $\vecw=(0,0)$, it is clear that $\vecz^*=(0,0)$. For $\vecw$ with $|\vecw|\neq 0$, $\vecz^*$ has the same direction as $\vecw$: $\vecz^{*}=\frac{|\vecz^{*}|}{|\vecw|}\vecw$. Thus to obtain $\vecz^{*}$, it is sufficient to calculate $|\vecz^{*}|$ as the minimizer of the following univariate problem:
\begin{equation}
\label{eq-admm-qsub-keyprob}
  \min\limits_{s \geq 0}\left\{\chi(s;\tau,\beta,t)=
  \min\big(\rho(s),\rho(\tau)\big)+\frac{\beta}{2}(s-t) ^2\right\},
\end{equation}
where $t=|\vecw|$.
In some special cases such as $\rho(s)=\rho_1(s)$, \eqref{eq-admm-qsub-keyprob} can be explicitly solved \cite{Fornasier2010Iterative}. For a general $\rho(\cdot)$ satisfying (AS1)(AS2)(AS3)(AS4), we derive an efficient computational framework.

For the convenience of description, we introduce the following two functions
\begin{align}
  \chi_1(s) = \rho(s)  + \frac{\beta}{2}(s - t)^2, \label{eq-admm-qsub-keyprob-1}\\
  \chi_2(s) = \rho(\tau) + \frac{\beta}{2}(s - t)^2. \label{eq-admm-qsub-keyprob-2}
\end{align}

\begin{proposition}
\label{proposition-admm-qsub-keyprob}
The minimization problem \eqref{eq-admm-qsub-keyprob} can be solved by
\begin{equation}
\label{eq-admm-qsub-keyprob-solution}
  s^* = \left\{
  \begin{array}{cll}
   s_1^*,  & \chi_1(s_1^*) < \chi_2(s_2^*),  \\
   \{s_1^*,s_2^*\},  & \chi_1(s_1^*) = \chi_2(s_2^*),  \\
   s_2^*,  & \chi_1(s_1^*) > \chi_2(s_2^*),
  \end{array}
  \right.
\end{equation}
where
$$
s_1^* = \argmin\limits_{0 \leq s \leq  \tau}\chi_1(s); s_2^* = \argmin\limits_{s \geq \tau}\chi_2(s) = \max\left(t, \tau \right).
$$
\end{proposition}
\proof It is clear that \eqref{eq-admm-qsub-keyprob-solution} gives the global minimizer of the problem \eqref{eq-admm-qsub-keyprob}. Because $\chi_2(s)$ is a quadratic univariate function, $s_2^* = \argmin\limits_{s \geq \tau}\chi_2(s) = \max\left(t, \tau \right)$.
\hfill$\blacksquare$

To compute $s_1^*$, we need the first and second order derivatives of $\chi_1$:
$$
\chi_1^{\prime}(s) = \rho^{\prime}(s) + \beta(s-t);\ \chi_1^{\prime\prime}(s) = \rho^{\prime\prime}(s) + \beta,
$$
and the following constant
$$
s_L = \inf\left\{s\in(0, + \infty): \rho^{\prime\prime}(s) > -\beta \right\}.
$$
Similarly to \cite{Nikolova2005Analysis,Chen2010Lower}, we have the following lower bound theory for the local minimizers of $\chi_1(s)$.
\begin{proposition}\label{proposition-secondorder-lowerbound}
[second order lower bound]
If $\rho(\cdot)$ satisfies (AS1)(AS2)(AS3)(AS4) and $s_{loc}^* $ is a local minimizer of $\min_{s \geq 0} \chi_1(s)$, then either $s_{loc}^*=0$ or $s_{loc}^* \geq s_L$.
\end{proposition}
\proof
The proof can be completed using the second order necessary condition, like \cite{Nikolova2005Analysis,Chen2010Lower}.
\hfill$\blacksquare$
\begin{remark}
\label{rem-lowerbound}
The constant $s_L$ (lower bound) for some regularization functions can be calculated explicitly.
For $\rho_1(s)= s$, we have $s_L = 0$. For $\rho_2(s)= s^p(0<p<1)$, we have $s_L = \left(\frac{p(1-p)}{\beta}\right)^{\frac{1}{2-p}}$. For $\rho_3(s) = \ln(\theta s + 1)$, we have $s_L = \max\left(0, \frac{1}{\sqrt{\beta}} - \frac{1}{\theta}\right)$. For $\rho_4(s) = \frac{\theta s}{\theta s + 1}$, we have $s_L = \max\left(0, \left(\frac{2}{\theta\beta}\right)^{\frac{1}{3}} - \frac{1}{\theta}\right)$.
\end{remark}
\begin{remark}
  When $s_L = 0$, the statement of Proposition~\ref{proposition-secondorder-lowerbound} is a degenerated result, which reads $s_{loc}^* = 0$ or $s_{loc}^* > 0$.
\end{remark}

Since $\rho(\cdot)$ satisfies (AS1)(AS2)(AS3)(AS4), $\chi_1(s)$ is concave in $[0, s_L]$ and strictly convex in $[s_L, +\infty)$. This observation is the key for us to design an efficient method to find the global minimizer of $\chi_1(s)$ over $[0,\tau]$.
\begin{proposition}\label{proposition-admm-qsub-keyprob-Chi1}
Under the assumptions of Proposition~\ref{proposition-secondorder-lowerbound}, we have:
\begin{itemize}
  \item[(1)] If $s_L >0, \chi_1^{\prime}(s_L) \geq 0$ ($s_L = 0, \chi_1^{\prime}(0+) \geq 0$), then $s_1^* = 0$ is the unique global minimizer of $\min_{0 \leq s \leq \tau} \chi_1(s)$.
  \item[(2)] If $s_L >0, \chi_1^{\prime}(s_L) < 0$ ($s_L = 0, \chi_1^{\prime}(0+) < 0$), then the equation $\chi_1^{\prime}(s) = 0$ has a unique root $\bar{s}$ on $[s_L,t]$. Set $\mathcal{X} = \{0,\min(\bar{s},\tau)\}$. The global minimizer of $\min_{0 \leq s \leq \tau} \chi_1(s)$ is given by $s_1^* = \arg\min_{s\in \mathcal{X}}\chi_1(s)$.
\end{itemize}
\end{proposition}
\proof
Two cases arise according to the value of $s_L$.
\begin{itemize}
  \item[(C1)] When $s_L = 0$, we have $\chi_1^{\prime\prime}(s) > 0$ on $(0,+\infty)$. Thus $\chi_1^{\prime}(s)$ strictly increases on $(0,+\infty)$.
  \begin{enumerate}
    \item If $\chi_1^{\prime}(0+) \geq 0$, we obtain $\chi_1^{\prime}(s) > 0$ on $(0,+\infty)$, that is $\chi_1(s)$ strictly increases on $(0,+\infty)$. Then $0$ is the unique global minimizer of $\min_{0 \leq s \leq \tau} \chi_1(s)$.
    \item If $\chi_1^{\prime}(0+) < 0$, then $\chi_1^{\prime}(s) = 0$ has a unique root $\bar{s}$ on $(0,+\infty)$ by the strict monotonicity of $\chi_1^{\prime}(s)$ and $\chi_1^{\prime}(s)\rightarrow +\infty $ as $s \rightarrow +\infty$. Actually, we have $0 = s_L \leq \bar{s} \leq t$, since for $\forall s > 0$, $\rho^{\prime}(s) \geq 0 $ and $\chi_1^{\prime}(\bar{s}) = \rho^{\prime}(\bar{s}) + \beta(\bar{s} - t) = 0$. Therefore, $\min(\bar{s},\tau)$ is the unique global minimizer of $\min_{0 \leq s \leq \tau} \chi_1(s)$.
  \end{enumerate}
  \item[(C2)] When $s_L > 0$, we have $\chi_1^{\prime\prime}(s) < 0$ on $(0,s_L)$ and $\chi_1^{\prime\prime}(s) > 0$ on $(s_L,+\infty)$. Thus $\chi_1^{\prime}(s)$ strictly decreases on $(0,s_L)$ and strictly increases on $(s_L,+\infty)$.
    \begin{enumerate}
      \item If $\chi_1^{\prime}(s_L) > 0$, then $\chi_1^{\prime}(s) > 0$ on $(0,+\infty)$, meaning that $\chi_1(s)$ strictly increases on $(0,+\infty)$. Then $0$ is the unique global minimizer of $\min_{0 \leq s \leq \tau} \chi_1(s)$.
      \item If $\chi_1^{\prime}(s_L) = 0$, then $s_L$ is a critical point but not a local minimizer. From $\chi_1^{\prime}(s) > 0$ on $(0,s_L)$ and $(s_L,+\infty)$, we can see that $\chi_1(s)$ strictly increases on $(0,+\infty)$. Then $0$ is the unique global minimizer of $\min_{0 \leq s \leq \tau} \chi_1(s)$.
      \item If $\chi_1^{\prime}(s_L) < 0$, then $\min_{s} \chi_1(s)$ has a unique local minimizer $\bar{s}$ on $[s_L,+\infty)$, i.e., $\chi_1^{\prime}(\bar{s}) = 0$. Since $\forall s > 0$, $\rho^{\prime}(s) \geq 0$, and $\chi_1^{\prime}(\bar{s}) = \rho^{\prime}(\bar{s}) + \beta(\bar{s} - t) = 0$, it follows that $s_L \leq \bar{s} \leq t$. Set $\mathcal{X} = \{0,\min(\bar{s},\tau)\}$. On the other hand, $\min_{s} \chi_1(s)$ has no local minimizer on $(0,s_L)$ due to Proposition \ref{proposition-secondorder-lowerbound}. Therefore, the global minimizer of $\min_{0 \leq s \leq \tau} \chi_1(s)$ is given by $\arg\min_{s\in \mathcal{X}}\chi_1(s)$.
    \end{enumerate}
\end{itemize}
Summarizing the conclusions of (C1) and (C2), we establish the proposition. 
\hfill$\blacksquare$

According to these discussions, we obtain an algorithm to find the global minimizer of \eqref{eq-admm-qsub-keyprob}; see \cref{alg-globalminimizer}.
\begin{algorithm}[t]
\caption{Find the global minimizer of the minimization problem \eqref{eq-admm-qsub-keyprob}}
\label{alg-globalminimizer}
\begin{algorithmic}[1]
    \REQUIRE $t$, $\tau$, the second order bound $s_L$ and functions $\chi_1(s)$, $\chi_1^{\prime}(s+)$, $\chi_2(s)$;
    \ENSURE $s^*$;
    \STATE $//$ Find the global minimizer of $s_1^* = \arg\min_{0 \leq s \leq  \tau} \chi_1(s)$.
    \IF{$\chi_1^{\prime}(s_L+) < 0$}
    \STATE Find the root $\bar{s}$ of equation $\chi_1^{\prime}(s) = 0$ in $[s_L,t]$;
    \STATE Set the feasible set $\mathcal{X} = \left\{0, \min(\bar{s}, \tau)\right\}$;
    \STATE Choose $s_1^* \in \mathcal{X}$ with $s_1^* := \arg\min_{s \in \mathcal{X}} \chi_1(s)$;
     \ELSE
    \STATE Set $s_1^* = 0$;
    \ENDIF
    \STATE $//$ Find the global minimizer of $s_2^* = \arg\min_{t \leq  \tau} \chi_2(s)$.
    \STATE Set $s_2^* = \max\{\tau,t\}$;
    \STATE $//$ Find the global minimizer $s^*$.
    \STATE Choose $s^*$ with
    $$
      s^* = \left\{
  \begin{array}{cll}
   s_1^*,  & \chi_1(s_1^*) < \chi_2(s_2^*),  \\
   \{s_1^*, s_2^*\},  & \chi_1(s_1^*) = \chi_2(s_2^*), \\
   s_2^*, & \mathrm{otherwise}.
  \end{array}
  \right.
    $$
\end{algorithmic}
\end{algorithm}

\begin{remark}\label{rem-root}
In \cref{alg-globalminimizer}, $\bar{s}\in [s_L,t]$ is the unique root of $\chi_1^{\prime}(s) = 0$ on the interval. It should be pointed out that the equation $\chi_1^{\prime}(s) = 0$ can be solved analytically in the cases $\rho_1(s)$ and $\rho_2(s)(p = \frac{1}{2},p = \frac{2}{3})$. In the more general case, no explicit formula has been found. Due to the uniqueness of the root, we can instead effectively numerically solve it in $[s_L,t]$ via any numerical procedure such as Brent's method~\cite{Brent1973Algorithm}.
\end{remark}

%

\subsection{Convergence Analysis}
The convergence of ADMM for non-convex composite problems is a difficult problem. Some results were obtained under some assumptions recently; see~\cite{Li2015Global,Wang2015Global} and the references therein. A key assumption is the surjectiveness of the linear mapping in the composite term, with a weakened form that the linear mappings in the constraint satisfies an image condition~\cite{Wang2015Global}. Unfortunately, these do not hold in our case, due to the non-surjectivess of the operator $\nabla=(\nabla_x,\nabla_y)$ and that a gradient vector field in 2D is required to satisfy the consistency condition. Their very important bounding techniques cannot be applied to our case.
Another very recent and important result is \cite{Guo2017Convergence} for the class of Douglas-Rachford splitting methods (DRSM) containing ADMM. Under mild assumptions, the authors proved the convergence of the DRSM for the ``strongly+weakly" convex programming. In Table~\ref{tab-PFs}, the regularizer functions $\rho_{1}$,  $\rho_{3}$, $\rho_{4}$ and $\rho_{8}$ are semi-convex. However, the truncated regularizer functions in general are not semi-convex. Currently we can only obtain the following result.

\begin{theorem} \label{theorem-2D-TR-ADMM-Convergence}
Assume that (AS1)(AS2)(AS3)(AS4)(AS5) hold and $\vecmu^{k+1}-\vecmu^{k}\rightarrow0$ as $k\rightarrow\infty$ in the ADMM in \cref{alg-admm}. Then any cluster point of the sequence $\{(\mathbf{u}^{k},\vecq^{k},\vecmu^{k})\}$, if exists, is a KKT point of the constrained optimization problem \cref{eq-2Diso-TR-equiv}.
\end{theorem}
\proof
First of all, we note that the regularization term $R(\cdot)$ is lower semi-continuous by (AS1)(AS2)(AS3)(AS4). This, together with (AS5), gives the properness, lower semi-continuity, and the coercivity of the augmented Lagrangian functional $\mathcal{L}(\mathbf{u},\vecq;\vecmu)$ for fixed Lagrange multiplier $\vecmu$.

By the ADMM, we have
\begin{equation}
\label{eq-admm-convergence-iterOptimalcondition}
\left\{
  \begin{aligned}
  &\partial R(\vecq^{k+1})+\vecmu^{k}+\beta(\vecq^{k+1}-\nabla\mathbf{u}^{k})\ni0,\\
  &\alpha \mathbf{A}^{\textsf{T}} (\mathbf{A} \mathbf{u}^{k+1}- \mathbf{f})-\nabla^{\textsf{T}}\vecmu^k - \beta\nabla^{\textsf{T}}\vecq^{k+1}- \beta\Delta \mathbf{u}^{k+1} =0,\\
  &\vecmu^{k+1} = \vecmu^{k}+\beta(\vecq^{k+1}- \nabla\mathbf{u}^{k+1}).
  \end{aligned}
\right.
\end{equation}
The assumption $\vecmu^{k+1}-\vecmu^{k}\rightarrow0$ as $k\rightarrow\infty$ gives immediately $\lim\limits_{k\rightarrow\infty}(\vecq^{k}- \nabla\mathbf{u}^{k})=0$. This assumption, together with the second equation in \cref{eq-admm-convergence-iterOptimalcondition} and the invertibility of $\mathbf{A}^{\textsf{T}}\mathbf{A}$, indicates that $\lim\limits_{k\rightarrow\infty}\mathbf{u}^{k+1}-\mathbf{u}^{k}=0$.
As $\vecq^{k+1}-\vecq^{k}=\vecq^{k+1}- \nabla\mathbf{u}^{k+1}+\nabla\mathbf{u}^{k+1}-\nabla\mathbf{u}^{k}+\nabla\mathbf{u}^{k}-\vecq^{k}$, we further deduce that $\lim\limits_{k\rightarrow\infty}\vecq^{k+1}-\vecq^{k}=0$.

Assume $(\mathbf{u}^{*},\vecq^{*},\vecmu^{*})$ to be a cluster point of $\{(\mathbf{u}^{k},\vecq^{k},\vecmu^{k})\}$. There is a subsequence $\{(\mathbf{u}^{k_i},\vecq^{k_i},\vecmu^{k_i})\}$ converging to $(\mathbf{u}^{*},\vecq^{*},\vecmu^{*})$. By the asymptotic regularity of $\{(\mathbf{u}^{k},\vecq^{k},\vecmu^{k})\}$, $\lim\limits_{i\rightarrow\infty}(\mathbf{u}^{k_i+1},\vecq^{k_i+1},\vecmu^{k_i+1})=(\mathbf{u}^{*},\vecq^{*},\vecmu^{*})$. Letting $k$ be $k_i$ in \cref{eq-admm-convergence-iterOptimalcondition} and sending $i$ to $\infty$ in the second and third equation, we get
\begin{equation}
\label{eq-admm-convergence-qiterOptimalcondition}
\partial R(\vecq^{k_i+1})+\vecmu^{k_i}+\beta(\vecq^{k_i+1}-\nabla\mathbf{u}^{k_i})\ni0,
\end{equation}
and
\begin{equation}
\label{eq-admm-convergence-limitofiterOptimalcondition}
\left\{
  \begin{aligned}
  &\alpha \mathbf{A}^{\textsf{T}} (\mathbf{A} \mathbf{u}^{*}- \mathbf{f})-\nabla^{\textsf{T}}\vecmu^*=0,\\
  &\vecq^{*}= \nabla\mathbf{u}^{*}.
  \end{aligned}
\right.
\end{equation}
As $\lim\limits_{i\rightarrow\infty}\vecmu^{k_i}=\vecmu^{*}$, $\lim\limits_{i\rightarrow\infty}\vecq^{k_i+1}=\vecq^{*}$, and $\lim\limits_{i\rightarrow\infty}(\vecq^{k_i+1}-\nabla\mathbf{u}^{k_i})=\lim\limits_{i\rightarrow\infty}(\vecq^{k_i+1}-\nabla\mathbf{u}^{k_i+1}+\nabla\mathbf{u}^{k_i+1}-\nabla\mathbf{u}^{k_i})=0$,
we have, by \cref{eq-subdifferential-closedness},
\begin{equation}
\label{eq-admm-convergence-limitofqiterOptimalcondition}
-\vecmu^*\in\partial R(\vecq^{*}),
\end{equation}
provided that $\lim\limits_{i\rightarrow\infty}R(\vecq^{k_i+1})=R(\vecq^{*})$.

We now show $\lim\limits_{i\rightarrow\infty}R(\vecq^{k_i+1})=R(\vecq^{*})$. On one hand, by letting $\mathcal{L}(\mathbf{u},\vecq;\vecmu) = \mathcal{L}(\mathbf{u}^{k_i+1},\vecq^{k_i+1};\vecmu^{k_i+1})$ in \cref{eq-2Diso-TR-equiv-augLagfun} and the
lower semi-continuity of $R$, as well as \cref{eq-admm-convergence-limitofiterOptimalcondition}, it holds that
\begin{equation}
\label{eq-admm-convergence-liminfR}
\liminf\limits_{i\rightarrow\infty}\mathcal{L}(\mathbf{u}^{k_i+1},\vecq^{k_i+1};\vecmu^{k_i+1})
\geq R(\vecq^{*})+\frac{\alpha}{2} \|\mathbf{A} \mathbf{u}^{*} - \mathbf{f}\|_V^2.
\end{equation}
On the other hand, in the ADMM, we have
\begin{equation}
\label{eq-admm-convergence-iterFunctionValue}
\left\{
  \begin{aligned}
  &\mathcal{L}(\mathbf{u}^{k},\vecq^{k+1};\vecmu^{k})\leq\mathcal{L}(\mathbf{u}^{k},\vecq^{*};\vecmu^{k}),\\
  &\mathcal{L}(\mathbf{u}^{k+1},\vecq^{k+1};\vecmu^{k})\leq\mathcal{L}(\mathbf{u}^{k},\vecq^{k+1};\vecmu^{k}),\\
  &\mathcal{L}(\mathbf{u}^{k+1},\vecq^{k+1};\vecmu^{k+1})=\mathcal{L}(\mathbf{u}^{k+1},\vecq^{k+1};\vecmu^{k})+\frac{1}{\beta}\|\vecmu^{k+1}-\vecmu^{k}\|_Q^2,
  \end{aligned}
\right.
\end{equation}
which yields
$$
\mathcal{L}(\mathbf{u}^{k+1},\vecq^{k+1};\vecmu^{k+1})\leq\mathcal{L}(\mathbf{u}^{k},\vecq^{*};\vecmu^{k})+\frac{1}{\beta}\|\vecmu^{k+1}-\vecmu^{k}\|_Q^2.
$$
Letting $k=k_i$ in the above and passing to limit, we obtain
\begin{equation}
\label{eq-admm-convergence-limsupR}
\limsup\limits_{i\rightarrow\infty}\mathcal{L}(\mathbf{u}^{k_i+1},\vecq^{k_i+1};\vecmu^{k_i+1})\leq R(\vecq^{*})+\frac{\alpha}{2} \|\mathbf{A} \mathbf{u}^{*} - \mathbf{f}\|_V^2.
\end{equation}
From \cref{eq-admm-convergence-liminfR} and \cref{eq-admm-convergence-limsupR}, it follows that $\lim\limits_{i\rightarrow\infty}R(\vecq^{k_i+1})=R(\vecq^{*})$. Then \cref{eq-admm-convergence-limitofqiterOptimalcondition} holds. \cref{eq-admm-convergence-limitofiterOptimalcondition} and \cref{eq-admm-convergence-limitofqiterOptimalcondition} shows that $(\mathbf{u}^{*},\vecq^{*},\vecmu^{*})$ is a KKT point of the constrained optimization problem \cref{eq-2Diso-TR-equiv}. ~\hfill$\blacksquare$

\begin{remark}\label{rem-2D-TR-DualVariableAsymptoticRegular}
For the assumption of $\vecmu^{k+1}-\vecmu^{k}\rightarrow0$ as $k\rightarrow\infty$ in \cref{theorem-2D-TR-ADMM-Convergence}, we have some numerical verification. Curves about the $\|\vecmu^{k+1}-\vecmu^{k}\|$ evolutions via $k$ of our experiments are shown in \cref{fig-residual}(c), which illustrate the reasonability of this assumption.
\end{remark}

\subsection{Numerical Experiments}



\graphicspath{{figures/}}

In this subsection we provide some numerical examples. We make comparisons between $\rho_1(s) = s$ (TV), $\rho_2(s) = s^p$ ($\ell_p, 0<p<1$), $\rho_3(s) = \ln(\theta s + 1)$ (LN) and $\rho_4(s) = \frac{\theta s}{\theta s + 1}$ (FRAC) and the corresponding truncated versions, TR-TV, TR-$\ell_p$, TR-LN and TR-FRAC, as well as  $\rho_{8}(s)$ (SCAD) and the truncated $\ell_2$ (TR-$\ell_2$), for the variational model \eqref{eq-generic-model-2Dani}. In all our numerical experiments, we choose $a = 3.7$ in the SCAD regularizer as suggested by Fan and Li ~\cite{Fan2001Variable}. The TR-$\ell_2$ is smooth at zero, while the others are nonsmooth at zero. Note we do not compare them to the standard $\ell_2$ regularization, because it is now well known that $\ell_2$ regularization cannot preserve image edges. The proposed ADMM is adopted to solve the minimization problems. Inspired by \cite{Dong2013Efficient} and lots of experiments, we use the following stopping condition
$$
\min\left\{\frac{||\bar{\mathbf{u}}^k -\bar{\mathbf{u}}^{k-1} ||_{V}}{||\mathbf{f}||_{V}},
\frac{||\bar{\mathbf{q}}^k -\nabla\bar{\mathbf{u}}^{k} ||_{Q}}{||\nabla\mathbf{f}||_{Q}} \right\} \leq 5 \times 10^{-5},
$$
where $\bar{\mathbf{u}}^k := \frac{1}{k+1} \sum_{j = 0}^{k} \mathbf{u}^j $ and $\bar{\mathbf{q}}^k := \frac{1}{k+1} \sum_{j = 0}^{k} \mathbf{q}^j $.

\begin{figure}[htbp]
  \captionsetup[subfigure]{justification=centering}
 \centering
 \newcommand{\imgwid}{0.25\textwidth}
  \begin{tabular}
  {c@{\hspace{3mm}}c@{\hspace{3mm}}c@{\hspace{3mm}}}
  \subfloat[][Shepp-Logan. Size:~$256\times 256$]
  {\includegraphics[width=\imgwid]{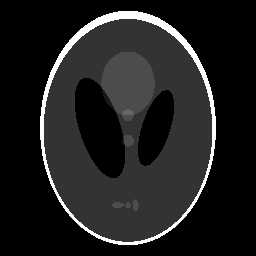}}&
  \subfloat[][TV. PSNR:~33.82dB ]
  {\includegraphics[width=\imgwid]{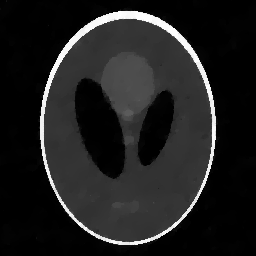}}&
    \subfloat[][SCAD. PSNR:~33.90dB ]
  {\includegraphics[width=\imgwid]{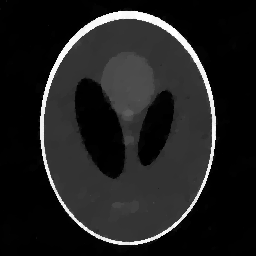}}
  \\
   \subfloat[][Noisy. PSNR:~20.20dB]
  {\includegraphics[width=\imgwid]{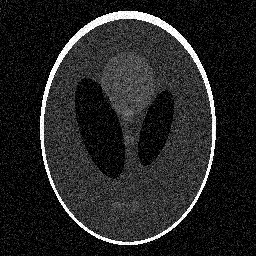}}&
  \subfloat[][TR-TV. PSNR:~36.99dB ]
  {\includegraphics[width=\imgwid]{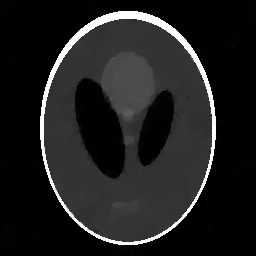}}&
   \subfloat[][TR-$\ell_2$. PSNR:~26.56dB ]
  {\includegraphics[width=\imgwid]{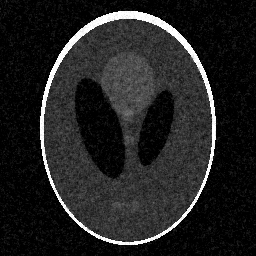}}
  \\
  \subfloat[][$\ell_p$. PSNR:~36.65dB]
    {\includegraphics[width=\imgwid]{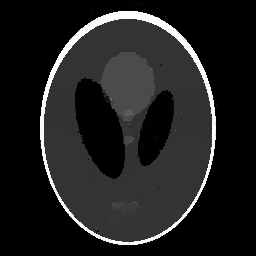}}&
  \subfloat[][LN. PSNR:~38.83B]
  {\includegraphics[width=\imgwid]{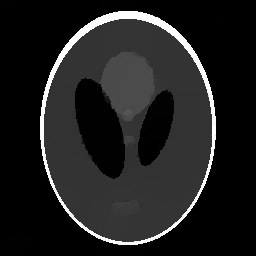}}&
  \subfloat[][FRAC. PSNR:~37.33dB ]
  {\includegraphics[width=\imgwid]{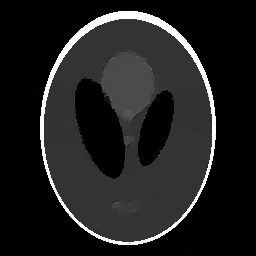}}
  \\
  \subfloat[][TR-$\ell_p$. PSNR:~36.93dB]
  {\includegraphics[width=\imgwid]{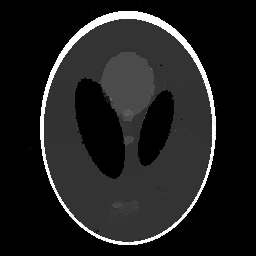}}&
  \subfloat[][TR-LN. PSNR:~39.18dB]
  {\includegraphics[width=\imgwid]{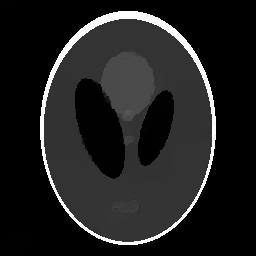}}&
  \subfloat[][TR-FRAC. PSNR:37.42dB]
  {\includegraphics[width=\imgwid]{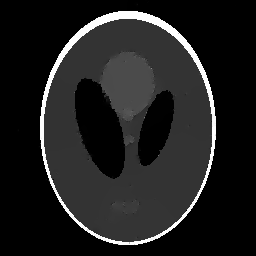}}
  \end{tabular}
  \caption{Image denoising ($\sigma = 25$) using different regularizer functions.
  The first and third rows: original \enquote{Shepp-Logan} image (a) and denoised images by (b) TV, (c) SCAD with $\theta = 1.0$, (g) $\ell_p$ with $p = 0.5$, (h) LN with $\theta = 10$ and (i) FRAC with $\theta = 10$; The second and fourth rows: noisy image  (d) and denoised images by (e) TR-TV with $\tau = 0.4$, (f) TR-$\ell_2$ with $\tau = 0.2$,  (j) TR-$\ell_p$ with $p = 0.5, \tau = 0.5$, (k) TR-LN with $\theta = 10, \tau = 0.5$, and (l) TR-FRAC with $\theta = 10, \tau = 0.5$. Note that $\ell_2$ regularizer over smoothes the edges, and therefore we do not put the result here.
  }
  \label{fig-Shepp-Logan-DN}
\end{figure}

\begin{figure}[htbp]
  \captionsetup[subfigure]{justification=centering}
 \centering
 \newcommand{\imgwid}{0.25\textwidth}
  \begin{tabular}
  {c@{\hspace{3mm}}c@{\hspace{3mm}}c@{\hspace{3mm}}}
  \subfloat[][QRcode. Size:~$378\times 378$]
  {\includegraphics[width=\imgwid]{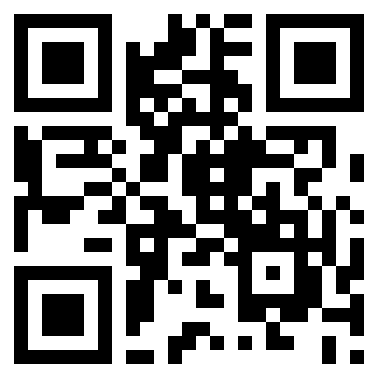}}&
  \subfloat[][TV. PSNR:~21.35dB ]
  {\includegraphics[width=\imgwid]{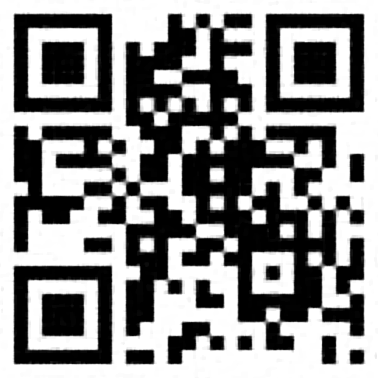}}&
    \subfloat[][SCAD. PSNR:~30.66dB ]
  {\includegraphics[width=\imgwid]{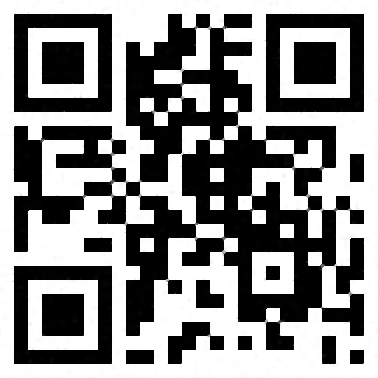}}
  \\
   \subfloat[][Blurry \& Noisy.  PSNR:~13.07dB]
  {\includegraphics[width=\imgwid]{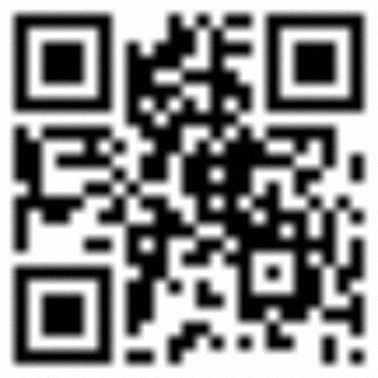}}&
  \subfloat[][TR-TV. PSNR: 29.94dB ]
  {\includegraphics[width=\imgwid]{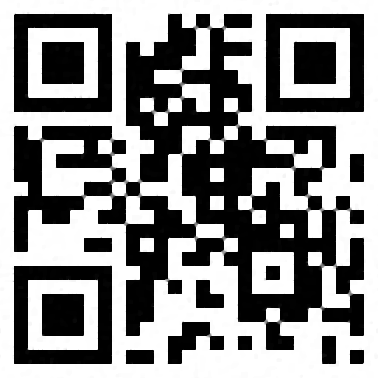}}&
   \subfloat[][TR-$\ell_2$.  PSNR:~27.20dB ]
  {\includegraphics[width=\imgwid]{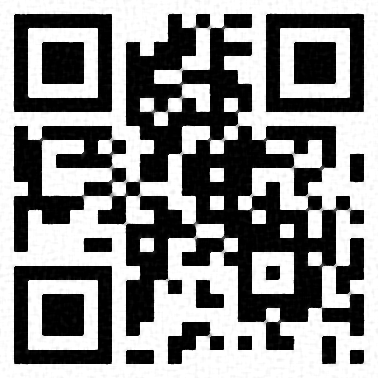}}
  \\
  \subfloat[][$\ell_p$. PSNR:~29.68dB]
    {\includegraphics[width=\imgwid]{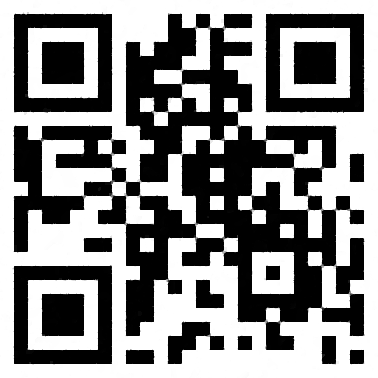}}&
  \subfloat[][LN. PSNR:~29.34dB]
  {\includegraphics[width=\imgwid]{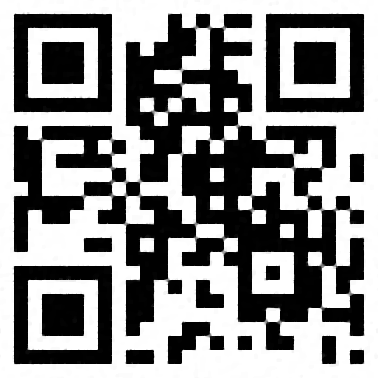}}&
  \subfloat[][FRAC. PSNR:~28.68dB ]
  {\includegraphics[width=\imgwid]{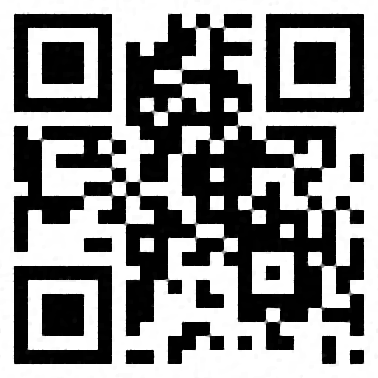}}
  \\
  \subfloat[][TR-$\ell_p$. PSNR:~30.87dB]
  {\includegraphics[width=\imgwid]{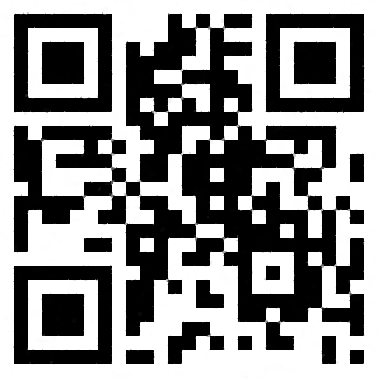}}&
  \subfloat[][TR-LN. PSNR:~30.56dB]
  {\includegraphics[width=\imgwid]{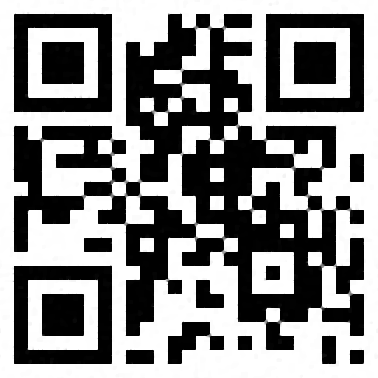}}&
  \subfloat[][TR-FRAC. PSNR:30.46dB]
  {\includegraphics[width=\imgwid]{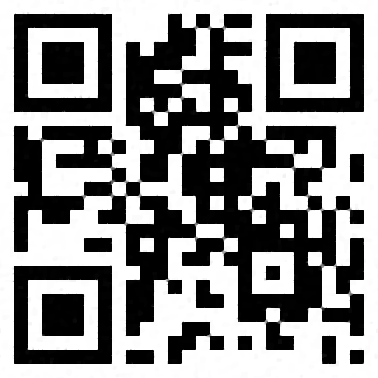}}
  \end{tabular}
  \caption{Image deblurring ($G(11,5), \sigma = 3$) using different regularizer functions.
  The first and third rows: original \enquote{QRcode} image (a) and deblurred images by (b) TV, (c) SCAD with $\theta = 0.2$, (g) $\ell_p$ with $p = 0.5$, (h) LN with $\theta = 1.0$ and (i) FRAC with $\theta = 1.0$; The second and fourth rows: blurry and noisy image  (d) and deblurred images by (e) TR-TV with $\tau = 0.5$, (f) TR-$\ell_2$ with $\tau = 0.2$, (j) TR-$\ell_p$ with $p = 0.5, \tau = 0.5$, (k) TR-LN with $\theta = 1.0, \tau = 0.5$, and (l) TR-FRAC with $\theta = 1.0, \tau = 0.5$. Note that $\ell_2$ regularizer over smoothes the edges, and therefore we do not put the result here.
  }
  \label{fig-QRcode-DB}
\end{figure}

\begin{figure}[htbp]
  \captionsetup[subfigure]{justification=centering}
 \centering
 \newcommand{\imgwid}{0.25\textwidth}
  \begin{tabular}
  {c@{\hspace{3mm}}c@{\hspace{3mm}}c@{\hspace{3mm}}}
  \subfloat[][Satellite. Size:~$135\times 135$]
  {\includegraphics[width=\imgwid]{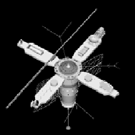}}&
  \subfloat[][TV. PSNR:~23.30dB ]
  {\includegraphics[width=\imgwid]{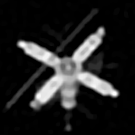}}&
    \subfloat[][SCAD. PSNR:~24.11dB ]
  {\includegraphics[width=\imgwid]{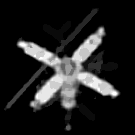}}
  \\
   \subfloat[][Blurry \& Noisy. PSNR:~19.99dB]
  {\includegraphics[width=\imgwid]{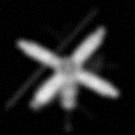}}&
  \subfloat[][TR-TV. PSNR:~23.95dB ]
  {\includegraphics[width=\imgwid]{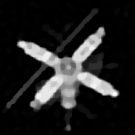}}&
   \subfloat[][TR-$\ell_2$.  PSNR:~23.01dB ]
  {\includegraphics[width=\imgwid]{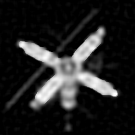}}
  \\
  \subfloat[][$\ell_p$. PSNR:~23.93dB]
    {\includegraphics[width=\imgwid]{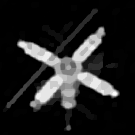}}&
  \subfloat[][LN. PSNR:~23.77dB]
  {\includegraphics[width=\imgwid]{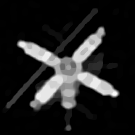}}&
  \subfloat[][FRAC. PSNR:~23.78dB ]
  {\includegraphics[width=\imgwid]{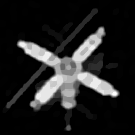}}
  \\
  \subfloat[][TR-$\ell_p$.  PSNR:~24.04dB]
  {\includegraphics[width=\imgwid]{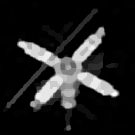}}&
  \subfloat[][TR-LN. PSNR:~24.04dB]
  {\includegraphics[width=\imgwid]{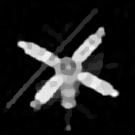}}&
  \subfloat[][TR-FRAC. PSNR:23.87dB]
  {\includegraphics[width=\imgwid]{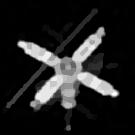}}
  \end{tabular}
  \caption{Image deblurring ($G(11,3), \sigma = 4$) using different regularizer functions.
  The first and third rows: original \enquote{Satellite} image (a) and deblurred images by (b) TV, (c) SCAD with $\theta = 0.1$, (g) $\ell_p$ with $p = 0.5$, (h) LN with $\theta = 10$ and (i) FRAC with $\theta = 10$; The second and fourth rows: blurry and noisy image  (d) and deblurred images by (e) TR-TV with $\tau = 0.5$, (f) TR-$\ell_2$ with $\tau = 0.2$, (j) TR-$\ell_p$ with $p = 0.5, \tau = 0.5$, (k) TR-LN with $\theta = 10, \tau = 0.5$, and (l) TR-FRAC with $\theta = 10, \tau = 0.5$. Note that $\ell_2$ regularizer over smoothes the edges, and therefore we do not put the result here.
  }
  \label{fig-Satellite-DB}
\end{figure}

We tested them using a standard \enquote{Shepp-Logan} image (\cref{fig-Shepp-Logan-DN}(a)), a \enquote{QRcode} image (\cref{fig-QRcode-DB}(a)) and a \enquote{Satellite} image (\cref{fig-Satellite-DB}(a)).  For image denoising problem, we did a test on \enquote{Shepp-Logan}; See \cref{fig-Shepp-Logan-DN}. The noisy observation is shown in \cref{fig-Shepp-Logan-DN} (d), which was corrupted by Gaussian noise with $\sigma = 25$.
The first and second rows show the denoising results by TV, SCAD, TR-TV and TR-$\ell_2$ regularizers. The third row shows the denoising results by $\ell_p (p=0.5)$, LN and FRAC regularizers; while the fourth row includes the results by TR-$\ell_p (p=0.5)$, TR-LN and TR-FRAC regularizers.
These results demonstrate that, for a regularization function, its truncated version can obtain better recovery with higher PSNR values. Besides, nonsmooth at zero functions generate much better denoising results than the smooth at zero TR-$\ell_2$.

\begin{figure}[htbp]
  \captionsetup[subfigure]{justification=centering}
 \centering
 \newcommand{\imgwid}{0.25\textwidth}
  \begin{tabular}
  {c@{\hspace{3mm}}c@{\hspace{3mm}}c@{\hspace{3mm}}}
  \subfloat[][Shepp-Logan. Size:~$256\times 256$]
  {\includegraphics[width=\imgwid]{Shepp-Logan-256.png}}&
  \subfloat[][TV. PSNR:~26.66dB ]
  {\includegraphics[width=\imgwid]{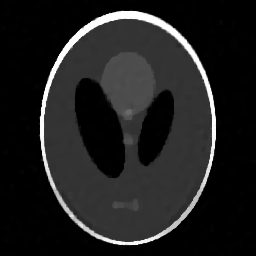}}&
    \subfloat[][SCAD. PSNR:~27.67dB ]
  {\includegraphics[width=\imgwid]{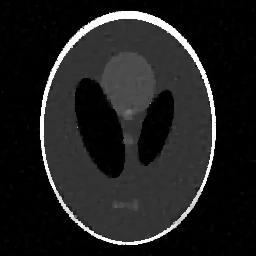}}
  \\
   \subfloat[][Blurry \& Noisy. PSNR:~19.02dB]
  {\includegraphics[width=\imgwid]{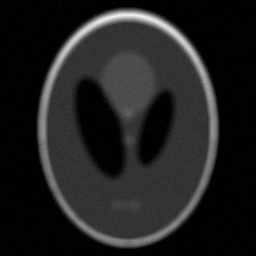}}&
  \subfloat[][TR-TV. PSNR:~27.38dB ]
  {\includegraphics[width=\imgwid]{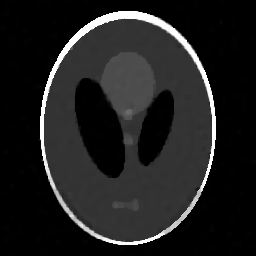}}&
      \subfloat[][TR-$\ell_2$. PSNR:~24.86dB ]
  {\includegraphics[width=\imgwid]{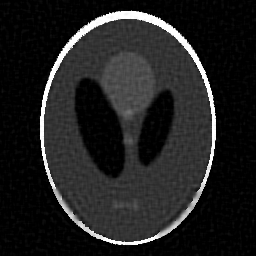}}
  \\
  \subfloat[][$\ell_p$. PSNR:~27.20dB]
    {\includegraphics[width=\imgwid]{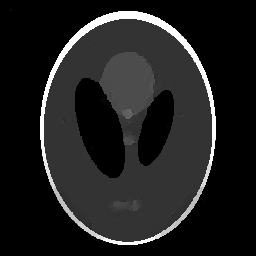}}&
  \subfloat[][LN. PSNR:~27.55dB]
  {\includegraphics[width=\imgwid]{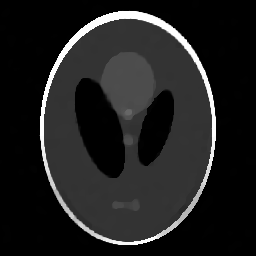}}&
  \subfloat[][FRAC. PSNR:~26.69dB ]
  {\includegraphics[width=\imgwid]{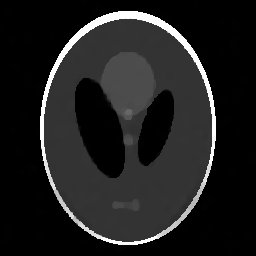}}
  \\
  \subfloat[][TR-$\ell_p$. PSNR:~27.45dB]
  {\includegraphics[width=\imgwid]{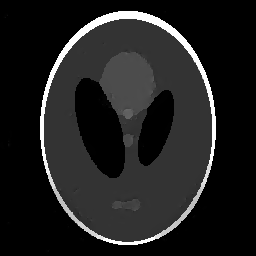}}&
  \subfloat[][TR-LN. PSNR:~28.02dB]
  {\includegraphics[width=\imgwid]{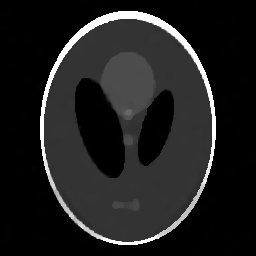}}&
  \subfloat[][TR-FRAC. PSNR:27.30dB]
  {\includegraphics[width=\imgwid]{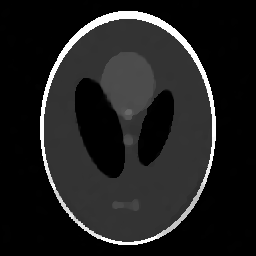}}
  \end{tabular}
  \caption{Image deblurring ($G(9,5), \sigma = 3$) using different regularizer functions.
  The first and third rows: original \enquote{Shepp-Logan} image (a) and deblurred images by (b) TV, (c) SCAD with $\theta = 0.1$, (g) $\ell_p$ with $p = 0.5$, (h) LN with $\theta = 1.0$ and (i) FRAC with $\theta = 1.0$; The second and fourth rows: blurry and noisy image  (d) and deblurred images by (e) TR-TV with $\tau = 0.7$, (f) TR-$\ell_2$ with $\tau = 0.2$, (j) TR-$\ell_p$ with $p = 0.5, \tau = 0.7$, (k) TR-LN with $\theta = 1.0, \tau = 0.7$, and (l) TR-FRAC with $\theta = 1.0, \tau = 0.7$. Note that $\ell_2$ regularizer over smoothes the edges, and therefore we do not put the result here.
  }
  \label{fig-Shepp-Logan-DB}
\end{figure}

For image deblurring problem, we tested our method on \enquote{Shepp-Logan}, \enquote{QRcode} and \enquote{Satellite}; See \cref{fig-Shepp-Logan-DB}, \cref{fig-QRcode-DB} and \cref{fig-Satellite-DB}. The observations were corrupted by different Gaussian blur and Gaussian noise, which are showed in~\cref{fig-Shepp-Logan-DB}(d) with $(G,9, 5)$ and $\sigma = 3$, \cref{fig-QRcode-DB}(d) with $(G,11, 5)$ and $\sigma = 3$ and~\cref{fig-Satellite-DB}(d) with $(G,11, 3)$ and $\sigma = 4$. In Figures~\ref{fig-Shepp-Logan-DB}, \ref{fig-QRcode-DB} and \ref{fig-Satellite-DB},
the first and second rows show the deblurring results by TV, SCAD, TR-TV and TR-$\ell_2$ regularizers. The third row shows the deblurring results by $\ell_p (p=0.5)$, LN and FRAC regularizers; while the fourth row includes the results by TR-$\ell_p (p=0.5)$, TR-LN and TR-FRAC regularizers.
These results demonstrate that, for a regularization function, its truncated version can obtain better recovery with higher PSNR values; See, especially the TR-TV.
Also, the SCAD regularizer performs very well.
Because there is some little noise in the observation, nonsmoothness at zero is still a benefit and TR-$\ell_2$ cannot outperform most nonsmooth at zero regularizers.
We also can see that all the truncated version models are more suitable for image deblurring problem.

\begin{figure}[htbp]
 \captionsetup[subfigure]{justification=centering}
  \centering
  \begin{tabular}{c@{\hspace{1mm}}c@{\hspace{1mm}}}
  \subfloat[][stopping condition $\frac{||\bar{\mathbf{q}}^k -\nabla\bar{\mathbf{u}}^{k} ||_{Q}}{||\nabla\mathbf{f}||_{Q}}$ ]
  {\includegraphics[width=0.45\textwidth]{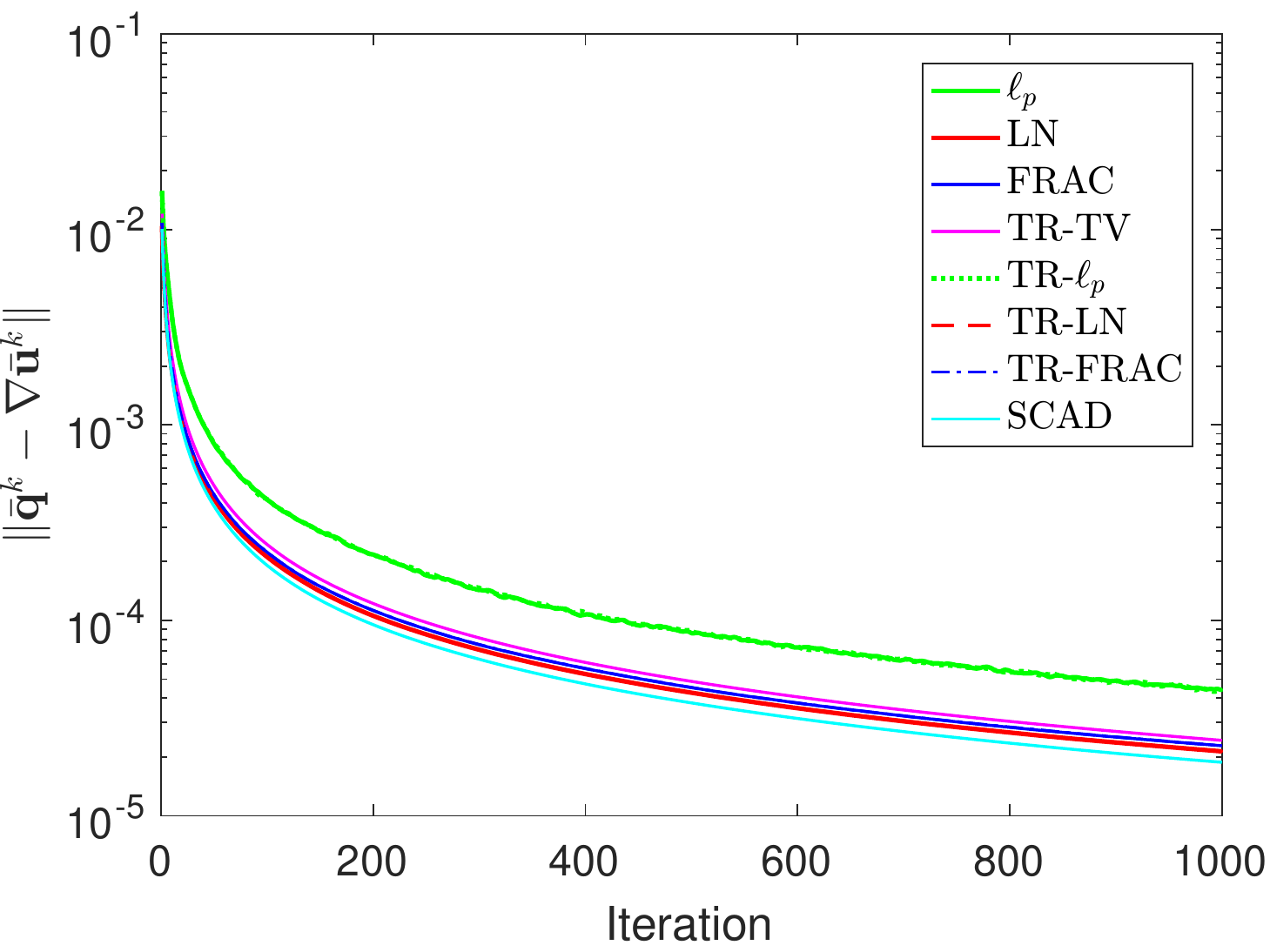}}&
  \subfloat[][stopping condition $\frac{||\bar{\mathbf{u}}^k -\bar{\mathbf{u}}^{k-1} ||_{V}}{||\mathbf{f}||_{V}}$ ]
  {\includegraphics[width=0.45\textwidth]{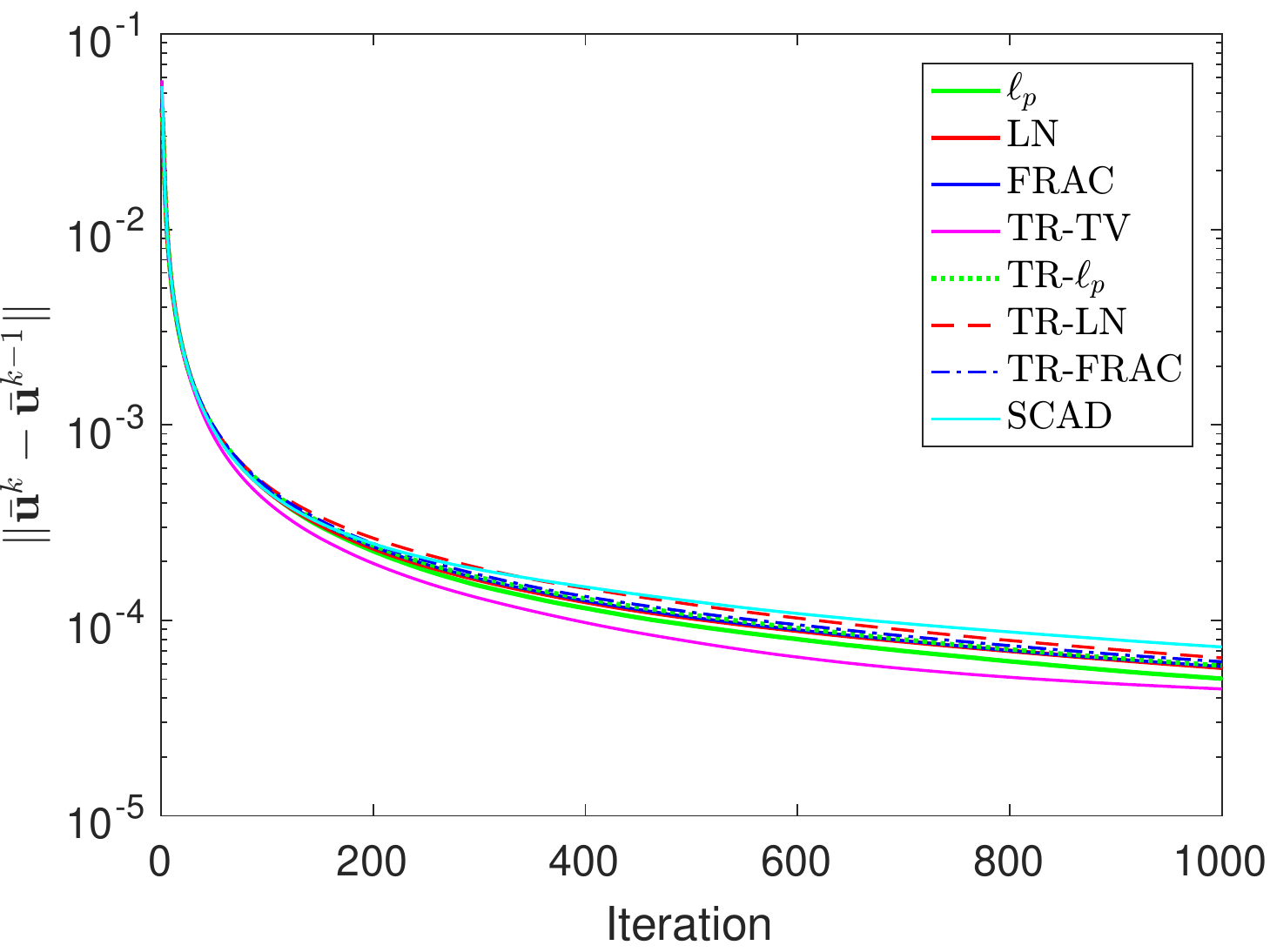}} \\
    \subfloat[][$\|\boldsymbol{\lambda}^{k+1} - \boldsymbol{\lambda}^k\|$ evolutions]
  {\includegraphics[width=0.45\textwidth]{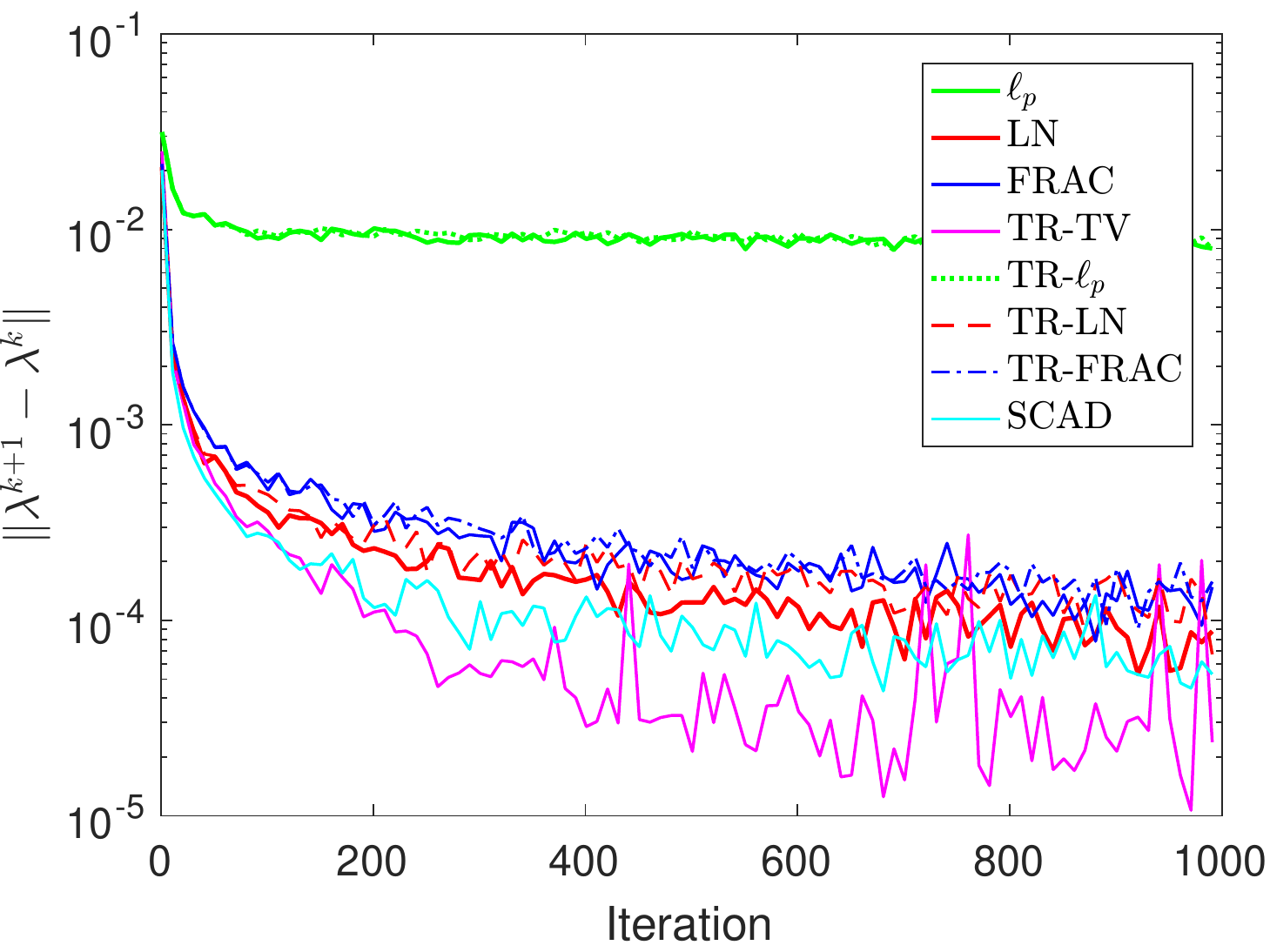}}&
   \end{tabular}
  \caption{The evolutions of  stopping condition and $\|\boldsymbol{\lambda}^{k+1} - \boldsymbol{\lambda}^k\|$  via iteration number $k$ for the \enquote{Satellite} image in \cref{fig-Satellite-DB}.}
  \label{fig-residual}
\end{figure}

We recorded the evolutions of stopping condition and $\|\boldsymbol{\lambda}^{k+1} - \boldsymbol{\lambda}^k\|$ via iteration number $k$ for the \enquote{Satellite} image in \cref{fig-residual}. We can see that ADMM is appropriate to solve the truncated regularization models.

Like most existing approaches, our method has parameters, including $\alpha,\beta$ and $\tau$ ($\theta$ for the SCAD regularizer). We did the experiments on almost all test images using different parameter values. We got the following observations. Firstly, $\alpha$ controls the restoration effect and smoothness of the result image. The smaller the $\alpha$ is, the smoother the restoration is. If $\alpha$ is extremely large (e.g., $1.e5$), the algorithm fails to restore the image. If $\alpha$ is extremely small (e.g., $1$), the image will become over-smoothed and some features will be lost. Actually, $\alpha$ is dependent on the noise level. Secondly, $\beta$ also controls the restoration effect. The larger $\beta$ is, the smoother the result is. Too large $\beta$ (e.g., $1.e5$) gives over-smoothed result. If $\beta$ is too small (e.g., $1$), there are some noise left in the result. Finally, $\tau$ ($\theta$ for the SCAD regularizer) relates to the lowest contrast of the images to preserve. To give analytical formulae of the optimal values for these parameters is quite difficult. Currently, we manually tuned these parameters according to the above observations, by starting from moderate values.
Due to the different shapes of regularizers, the optimal parameters of our method with different regularizers are different.

In Table~\ref{tab-Param-alpha}, Table~\ref{tab-Param-beta} and Table~\ref{tab-Param-tau}, we recorded the PSNR values for different $\alpha$, $\beta$ and $\tau$ ($\theta$ for the SCAD regularizer) with fixed other parameters, respectively. We can see that for each regularizer, there are certain reasonably large intervals for its parameters generating good restoration results. Of course, different regularizers may have different effective parameter intervals, which is reasonable, due to, again, the different shapes of regularizer functions.

\begin{table}[htbp]
  \centering
  \caption{The PSNR values under different $\alpha$ in the proposed ADMM for the \enquote{Satellite} image.}
    \begin{tabular}{|cc|cc|cc|cc|cc|}
    \hline
    \multicolumn{2}{|c|}{TR-TV}  & \multicolumn{2}{|c|}{TR-$\ell_p$} & \multicolumn{2}{|c|}{TR-LN} & \multicolumn{2}{|c|}{TR-FRAC} & \multicolumn{2}{|c|}{SCAD}\\
    \hline
    \multicolumn{2}{|c|}{$(\beta,\tau)=$}  & \multicolumn{2}{|c|}{$(\beta,p,\tau)=$} & \multicolumn{2}{|c|}{$(\beta,\theta,\tau)=$} & \multicolumn{2}{|c|}{$(\beta,\theta,\tau)=$} & \multicolumn{2}{|c|}{$(\beta,\theta)=$}\\
    \multicolumn{2}{|c|}{$(600,0.6)$}  & \multicolumn{2}{|c|}{$(5000,0.5,0.2)$} &
    \multicolumn{2}{|c|}{$(8000,10,0.2)$} & \multicolumn{2}{|c|}{$(6000,10,0.2)$} &
    \multicolumn{2}{|c|}{$(100,0.1)$}\\
     \hline
     $\alpha$ & PSNR & $\alpha$ & PSNR & $\alpha$ & PSNR & $\alpha$ & PSNR & $\alpha$ & PSNR \\
    \hline
    1200  & 23.91  & 1000  & 23.43  & 9000  & 23.99  & 4000  & 23.70  & 50    & 23.67  \\
    1400  & 23.95  & 2000  & 23.86  & 10000  & 24.01  & 5000  & 23.78  & 100   & 23.97  \\
    1600  & 23.98  & 3000  & 23.98  & 11000  & 24.03  & 6000  & 23.83  & 150   & 24.08  \\
    1800  & 24.00  & 4000  & 24.08  & 12000  & 24.04  & 7000  & 23.86  & 200   & 24.11  \\
    2000  & 24.00  & 5000  & 24.04  & 13000  & 24.04  & 8000  & 23.87  & 250   & 24.07  \\
    2200  & 24.00  & 6000  & 24.02  & 14000  & 24.04  & 9000  & 23.87  & 300   & 23.98  \\
    2400  & 23.98  & 7000  & 23.95  & 15000  & 24.04  & 10000  & 23.86  & 350   & 23.86  \\
    2600  & 23.95  & 8000  & 23.87  & 16000  & 24.02  & 11000  & 23.84  & 400   & 23.70  \\
    2800  & 23.92  & 9000  & 23.77  & 17000  & 24.01  & 12000  & 23.82  & 450   & 23.50  \\
    3000  & 23.87  & 10000  & 23.63  & 18000  & 23.99  & 13000  & 23.77  & 500   & 23.27  \\
    \hline
    \end{tabular}%
  \label{tab-Param-alpha}%
\end{table}%

\begin{table}[htbp]
  \centering
  \caption{The PSNR values under different $\beta$ in the proposed ADMM for the \enquote{Satellite} image.}
    \begin{tabular}{|cc|cc|cc|cc|cc|}
    \hline
    \multicolumn{2}{|c|}{TR-TV}  & \multicolumn{2}{|c|}{TR-$\ell_p$} & \multicolumn{2}{|c|}{TR-LN} & \multicolumn{2}{|c|}{TR-FRAC} & \multicolumn{2}{|c|}{SCAD}\\
    \hline
    \multicolumn{2}{|c|}{$(\alpha,\tau)=$} & \multicolumn{2}{|c|}{$(\alpha,p,\tau)=$} & \multicolumn{2}{|c|}{$(\alpha,\theta,\tau)=$} & \multicolumn{2}{|c|}{$(\alpha,\theta,\tau)=$} & \multicolumn{2}{|c|}{$(\alpha,\theta)=$}\\
    \multicolumn{2}{|c|}{$(2000,0.6)$} & \multicolumn{2}{|c|}{$(5000,0.5,0.2)$} & \multicolumn{2}{|c|}{$(13000,10,0.2)$} & \multicolumn{2}{|c|}{$(8000,10,0.2)$}  & \multicolumn{2}{|c|}{$(200,0.1)$}\\
    \hline
    $\beta$  & PSNR & $\beta$  & PSNR & $\beta$  & PSNR & $\beta$  & PSNR & $\beta$  & PSNR \\
    \hline
    200   & 23.49  & 1000  & 22.96  & 4000  & 23.43  & 2000  & 22.95  & 20    & 22.02  \\
    300   & 23.71  & 2000  & 23.54  & 5000  & 23.66  & 3000  & 23.28  & 40    & 22.73  \\
    400   & 23.93  & 3000  & 23.82  & 6000  & 23.88  & 4000  & 23.63  & 60    & 23.56  \\
    500   & 23.99  & 4000  & 23.95  & 7000  & 24.01  & 5000  & 23.83  & 80    & 24.01  \\
    600   & 24.00  & 5000  & 24.04  & 8000  & 24.04  & 6000  & 23.87  & 100   & 24.11  \\
    700   & 23.99  & 6000  & 24.04  & 9000  & 24.01  & 7000  & 23.85  & 120   & 24.01  \\
    800   & 23.95  & 7000  & 23.94  & 10000  & 23.95  & 8000  & 23.78  & 140   & 23.85  \\
    900   & 23.91  & 8000  & 23.82  & 11000  & 23.86  & 9000  & 23.68  & 160   & 23.70  \\
    1000  & 23.86  & 9000  & 23.73  & 12000  & 23.77  & 10000  & 23.58  & 180   & 23.55  \\
    1100  & 23.80  & 10000  & 23.63  & 13000  & 23.68  & 11000  & 23.48  & 200   & 23.42  \\
    \hline
    \end{tabular}%
  \label{tab-Param-beta}%
\end{table}%

\begin{table}[htbp]
  \centering
  \caption{The PSNR values under different truncation parameter $\tau$ ($\theta$ for the SCAD regularizer) in the proposed ADMM for the \enquote{Satellite} image. Note that, in our experiments, the intensity of image data was normalized to be $[0,1]$. For the SCAD regularizer, we chose $a = 3.7$ as suggested by Fan and Li~\cite{Fan2001Variable}.}
    \begin{tabular}{|cc|cc|cc|cc|cc|}
    \hline
    \multicolumn{2}{|c|}{TR-TV} & \multicolumn{2}{|c|}{TR-$\ell_p$} & \multicolumn{2}{|c|}{TR-LN} & \multicolumn{2}{|c|}{TR-FRAC} & \multicolumn{2}{|c|}{SCAD}  \\
    \hline
    \multicolumn{2}{|c|}{$(\alpha,\beta)=$} & \multicolumn{2}{|c|}{$(\alpha,\beta,p)=$} & \multicolumn{2}{|c|}{$(\alpha,\beta,\theta)=$} & \multicolumn{2}{|c|}{$(\alpha,\beta,\theta)=$} & \multicolumn{2}{|c|}{$(\alpha,\beta)=$} \\
    \multicolumn{2}{|c|}{$(2000,600)$} & \multicolumn{2}{|c|}{$(5000,5000,0.5)$} &
    \multicolumn{2}{|c|}{$(13000,8000,10)$} & \multicolumn{2}{|c|}{$(8000,6000,10)$} &
     \multicolumn{2}{|c|}{$(200,100)$} \\
    \hline
    $\tau$  & PSNR & $\tau$  & PSNR & $\tau$  & PSNR & $\tau$  & PSNR & $\theta(\tau = a \theta)$  & PSNR \\
    \hline
    0.1   & 23.12  & 0.1   & 23.63  & 0.1   & 23.57  & 0.1   & 23.52  & 0.01  & 22.39  \\
    0.2   & 23.58  & 0.2   & 24.04  & 0.2   & 24.04  & 0.2   & 23.87  & 0.04  & 23.22  \\
    0.3   & 23.48  & 0.3   & 23.98  & 0.3   & 23.90  & 0.3   & 23.83  & 0.07  & 23.89  \\
    0.4   & 23.81  & 0.4   & 23.93  & 0.4   & 23.90  & 0.4   & 23.81  & 0.10  & 24.11  \\
    0.5   & 23.99  & 0.5   & 23.93  & 0.5   & 23.90  & 0.5   & 23.81  & 0.13  & 24.07  \\
    0.6   & 24.00  & 0.6   & 23.93  & 0.6   & 23.90  & 0.6   & 23.81  & 0.16  & 23.93  \\
    0.7   & 24.00  & 0.7   & 23.93  & 0.7   & 23.90  & 0.7   & 23.81  & 0.19  & 23.77  \\
    0.8   & 24.00  & 0.8   & 23.93  & 0.8   & 23.90  & 0.8   & 23.81  & 0.22  & 23.61  \\
    0.9   & 24.00  & 0.9   & 23.93  & 0.9   & 23.90  & 0.9   & 23.81  & 0.25  & 23.49  \\
    1.0   & 24.00  & 1.0   & 23.93  & 1.0   & 23.90  & 1.0   & 23.81  & 0.27  & 23.48  \\
 \hline
 \end{tabular}%
  \label{tab-Param-tau}%
\end{table}%

\section{Conclusions}
Signal and image restoration is a typical inverse problem and energy regularization methods have been proven very useful. Various regularization functions have been proposed in the literature, including convex ones, smooth ones, and nonconvex or nonsmooth ones. In this paper, we studied the possibility of contrast-preserving restoration by variational methods. It is shown that any convex or smooth regularization, or even most nonconvex regularization, is impossible or with very low probability to recover the ground truth. We then naturally presented a general regularization framework based on truncation. Some analysis in 1D theoretically demonstrate its better contrast-preserving ability. Optimization algorithms in 2D with implementation details and convergence verification were given. Experiments numerically showed the advantages of the framework. One future work is to design more efficient algorithms to solve the models.

\Acknowledgements{The authors thank Dr. Chao Zeng to proof read the paper.
This work was supported by National Natural Science Foundation of China (Grant No. 11301289 and 11531013).}


\end{document}